\documentclass[11pt]{article}
\usepackage{amsfonts,amssymb,amsmath,amsgen}
\usepackage{latexsym}
\newcommand{\Hi}{{\cal H}}
\newcommand{\K}{{\cal K}}
\newcommand{\RR}{{\mathbb R}}
\newcommand{\ZZ}{{\mathbb Z}}
\newcommand{\Ft}{\mathcal{F}_t}
\begin{document}
\newtheorem{theorem}{THEOREM}{\bf}{\it}
\newtheorem{lemma}{LEMMA}[section]{\bf}{\it}
\newtheorem{prop}{PROPOSITION}{\bf}{\it}
\newtheorem{corol}{COROLLARY}{\bf}{\it}
\newtheorem{remark}{Remark}{\bf}{\rm}
\newcommand{\cqfd}{\mbox{ } \hfill$\Box$}

\title{
Strichartz estimates\\ for the Wave and Schr\"odinger  Equations\\
with the Inverse-Square Potential}

\author{
Nicolas Burq$^{\mathbf a}$
\and
Fabrice Planchon$^{\mathbf b}$
\and
John G. Stalker$^{\mathbf c}$
\and
A. Shadi Tahvildar-Zadeh$^{\mathbf d}$
}

\date{\today}

\maketitle
\begin{abstract}
We prove spacetime weighted-$L^2$ estimates for the Schr\"odinger and wave
equation with an inverse-square potential.
We then deduce Strichartz estimates for these equations.
\end{abstract}
\section{Introduction}

Consider the following linear equations
\begin{eqnarray}\label{eq:lse}
i\partial_t u +\Delta u- \frac{a}{|x|^2} u = 0, &\qquad &
u(0,x) = f(x)
\\
\label{eq:lwe}
-\partial_t^2 u +\Delta u- \frac{a}{|x|^2} u = 0, &\qquad &
u(0,x) = f(x),\ \partial_t u(0,x) = g(x)
\end{eqnarray}
where $\Delta $ is the $n$ dimensional Laplacian and $a$ is a real
number. The Schr\"odinger equation (\ref{eq:lse}) is of interest in
quantum mechanics (see \cite{Cas50,KalSchWalWus75,CEFG01} and references therein),
while the wave equation (\ref{eq:lwe}) arises in the study of  wave
propagation on conic manifolds \cite{CheTay82a}. We also note that the
heat flow for the elliptic operator $-\Delta + a|x|^{-2}$ has been studied
in the theory of combustion (see \cite{VazZua00} and references
therein). 

The mathematical interest in these equations however comes mainly from
the fact that the potential term is homogeneous of degree -2 and
therefore scales exactly the same as the Laplacian. This in particular
implies that perturbation methods cannot be used in studying the effect
of this potential. Indeed, the $|x|^{-2}$ decay is in some sense the borderline case
for the existence of global-in-time estimates for wave or Schr\"odinger
equations with a potential \cite{RodS,GeoVis}. In particular, it is known that a {\em negative} potential $V$ decaying
slower than inverse-square results in the spectrum of $-\Delta+V$ being
unbounded from below \cite[\S XIII, pp.87--88]{ReedSimon3}.

On the other hand, the scale-covariance present in the elliptic operator $-\Delta+\frac{a}{|x|^2}$ appearing in the
 above equations is a key
feature of many problems in physics and in geometry, where such scaling
behavior manifests itself if not everywhere, at least in a certain
region of space, for example near a singular point, or near infinity, or
both. Consider for example the Dirac equation with a Coulomb potential,
(which can be recast in the form of a Klein-Gordon equation with an
inverse-square potential, plus other terms which exhibit the same
scaling as the Laplacian) \cite{Cas50}. Another family of examples is given by linearized perturbations
of spacetime metrics that
are well-known solutions of the Einstein equations of general
relativity, such as the Schwarzschild solution \cite{RegWhe57,Zer70} or the
Reissner-Nordstr\"om solution \cite{Zer74, Mon74}.

 There are also  {\em nonlinear} problems, of a geometric nature, where such critical behavior potentials make an
appearance, for example the perturbation of equivariant stationary
solutions of a 2+1-dimensional wave map from the Minkowski space into a
2-sphere gives rise to a system of two wave equations, with a potential
that behaves like $|x|^{-2}$ both near zero and near infinity, as well
as another term with the same scaling, i.e $|x|^{-2}x\cdot
\nabla$. The occurrence of this phenomenon in a nonlinear setting is
significant, since it is clear that to study a nonlinear wave equation
one must have estimates for the linear inhomogeneous wave equation that bound various (perhaps fractional)
number of derivatives of the solution in term of the correct number of
derivatives of the source and the  data. Such estimates are by-and-large unknown
for problems involving a potential, except for  those that are of much faster decay than
$|x|^{-2}$ (\cite{Yaj95b} and references therein). 

Equations such as (\ref{eq:lse}) and (\ref{eq:lwe}) with the inverse-square
potential thus represents the simplest case, where the scaling holds
exactly and everywhere. These are to be thought of as model problems,
used to develop and test new tools and methods that we hope are capable of
being generalized to the more complicated situations that are of actual
physical and geometric interest, such as those named above. 

In \cite{PST1} we showed for the wave equation (\ref{eq:lwe}) that in
the radial case, i.e. when the data -- and thus the solution -- are
radially symmetric, the solution to (\ref{eq:lwe}) satisfies generalized
spacetime Strichartz estimates as long as
\begin{equation}\label{cond:a}
a> -(n-2)^2/4.
\end{equation}
The corresponding Strichartz estimates would hold for the Schr\"odinger
equation (\ref{eq:lse}) as well, since our proof was based solely on estimates for
the elliptic operator
\[
P_a := -\Delta+\frac{a}{|x|^2}.
\]
\begin{remark} {\rm  $P_a$ is in fact  the {\em self-adjoint extension} of $-\Delta+a|x|^{-2}$.
It is known that in the range $-(n-2)^2/4<a<1-(n-2)^2/4$ the extension is not unique
\cite{Tit46,KalSchWalWus75}. In this case, when we do make a choice among the possible 
extensions, such as in (\ref{nupos}), it corresponds to the Friedrichs extension
(see \cite{PST1,KalSchWalWus75} for details).}
\end{remark}

In this work we intend to remove the assumption on the data being radially
symmetric. As explained in \cite{PST1,PST2}, one cannot hope to get any
kind of dispersive (be it at fixed $t$ or spacetime) estimate if
(\ref{cond:a}) is not satisfied. We also note that when $a<0$, the
classical $L^\infty-L^1$ estimate for the wave equation does not hold
\cite{PST2}, and thus one cannot obtain Strichartz estimates by
interpolation between this dispersive estimate and the energy estimate
(see also the remark at the end of Subsection \ref{StriSchrod}).

This paper is divided into four sections. In Section~2, we obtain
weighted-$L^2$ estimates for (\ref{eq:lse}) and (\ref{eq:lwe}). Such
estimates are known for the free Schr\"odinger equation, and are often
referred to as local smoothing estimates (\cite{BAK} and references
therein). In Section~3 we deduce Strichartz estimates for solutions of
the Schr\"odinger equation (\ref{eq:lse}) through Duhamel's formula,
combining the smoothing estimate (\ref{eq:smoothS}) for (\ref{eq:lse})
with Strichartz estimates for the free Schr\"odinger
equation. We then do the same for the wave equation. We note that for the Schr\"odinger equation 
such a strategy was successfully applied in
\cite{JouSofSog91} for rapidly decaying potentials, to obtain dispersive
estimates (from which Strichartz estimates are then obtained by the
usual duality argument). This result was recently extended to potentials
that decay strictly faster than $|x|^{-2}$ \cite{RodS}. While our potential
obviously does not satisfy this condition, we are able to take advantage
of its special form to extend the approach in \cite{RodS} to our
setting.  Finally, in
Section~4 we obtain the frequency-localized version of the above
estimates and use them to obtain the generalized Strichartz estimates
(with derivatives) for these equations, which we then apply to obtain an optimal
global well-posedness result for a nonlinear wave equation.
\subsection{Notations}
In this paper we will be using the following notations.  For integer $n\geq 2$ Let 
\begin{equation}\label{def:lambda}
\lambda(n) := \frac{n-2}{2}.
\end{equation}
For integer $d\geq 0$ and real number $a\geq -\lambda^2(n)$ let
\begin{equation}\label{def:nud}
\nu_d(n,a) = \sqrt{(\lambda(n)+d)^2+a}
\end{equation}
We will suppress the arguments of the above functions whenever doing so does not cause confusion.

We  also define multiplication operators $\Omega^s$ by
\begin{displaymath}
(\Omega^s \phi)(x) = |x|^s \phi (x).
\end{displaymath}
Abusing notation we use the same symbol for the operators which are pointwise
equal to these for all times,
\begin{displaymath}
(\Omega^s \phi)(t,x) = |x|^s \phi (t,x).
\end{displaymath}
For $|s|<n/2$ and integer $d\geq 0$ let $\dot{H}^s_{\geq d}$ denote the subspace of the homogeneous 
Sobolev 
space $\dot{H}^s(\RR^n)$ consisting of functions that are orthogonal to all spherical harmonics of degree 
less than $d$, and let $\dot{H}^s_{<d}$ denote the orthogonal complement of this space.  Finally, let
\[
d_0(n) := \left\{\begin{array}{ll} 1 & n=2 \\ 0 & n\geq 3.\end{array}\right.
\]
\section{Weighted-$L^2$ estimates for the Schr\"odinger and wave equations}
\subsection{Local smoothing for the Schr\"odinger equation}

Except for the case $n=2$ the following theorem is well-known when $a=0$.  See
\cite{Sim92} for the sharp constants and the references therein for the
history.  In this paper we need only the case $\alpha=1/4$.
\begin{theorem}\label{thm:smoothing}
Let $n\geq 2$, $d\geq d_0(n)$, $0<\alpha<\frac{1}{4} + \frac{1}{2}\nu_d$ and let $u$ be the unique 
solution of {\rm(\ref{eq:lse})}.  
There exists a constant $C>0$, depending on
$n$, $a$, $d$ and $\alpha$, such that for all $f\in L^2_{\geq d}(\RR^n)$,
\begin{equation}\label{eq:smoothS}
\| \Omega^{-1/2-2\alpha}P_a^{1/4-\alpha} u \|_{L^2(\RR^{n+1})}
\le C \| f \|_{L^2 (\RR^n)}
\end{equation}
\end{theorem}
{\em Proof of Theorem~\ref{thm:smoothing}}:
We begin by noting that, by rotational symmetry and the $L^2$ orthogonality
of the various spherical harmonic spaces, it suffices to prove the
estimate~(\ref{eq:smoothS}) for $f$ belonging to the $l$'th harmonic subspace,
where $l\ge 0$ if $\lambda^2 + a > 0$ and $l>0$ if $\lambda^2 + a = 0$.
This, of course, requires the constants to be uniformly bounded in $l$, but
this will be clear from the explicit form of the constants given below.
On the $l$'th spherical harmonic subspace
\begin{equation}
P_a = A_\nu,
\end{equation}
where
\begin{eqnarray}
A_\nu &=& -\partial^2_r -(n-1)r^{-1}\partial_r + [l(l+2\lambda) + a]r^{-2}\nonumber\cr
      &=& -\partial^2_r - (n-1)r^{-1}\partial_r + [\nu^2 - \lambda^2]r^{-2},\label{def:Anu}
\end{eqnarray}
where
\begin{equation}
\nu = \nu_l(n,a)= \sqrt{(\lambda + l)^2 + a}.
\end{equation}
Our assumptions imply that
\begin{equation}\label{nupos}
\nu > 0.
\end{equation}

The above considerations allow us to restate the problem as follows.
We are to prove the estimate
\begin{displaymath}
\|\Omega^{-1/2-2\alpha}A^{1/4-\alpha}_\nu S_\nu f\|_{L^2(\RR^{n+1})} \le C
\|f\|_{L^2(\RR^n)}
\end{displaymath}
where $S_\nu f$ is the unique solution of the initial value problem
\begin{displaymath}
i\partial_t u - A_\nu u = 0, \qquad u(0,x) = f(x).
\end{displaymath}
In this we are allowed to assume that $f$ belongs to the $l$'th spherical
harmonic subspace, but in fact we have no further use for this assumption.

We define the Hankel transform of order $\nu$
in the usual way:
\begin{displaymath}
(\Hi_\nu \phi)(\xi) = \int^\infty_0 (r|\xi|)^{-\lambda}
J_\nu(r|\xi|) \phi(r\xi/|\xi|) r^{n-1}\,dr,
\end{displaymath}
where $J_\nu$ is the Bessel function of the first kind of order $\nu$.  
By abuse of notation we use the same symbol $\Hi_\nu$ to denote the operator on
functions on $\RR^{n+1}$ which is just the Hankel transform of order~$\nu$
pointwise for all times:
\begin{displaymath}
(\Hi_\nu \phi)(t,\xi) = \int^\infty_0 (r|\xi|)^{-\lambda}
J_\nu(r|\xi|) \phi(t,r\xi/|\xi|) r^{n-1}\,dr.
\end{displaymath}
The Hankel transform has the following properties:
\begin{itemize}
\item[(i)] $\Hi_\nu^2 = 1$,
\item[(ii)] $\Hi_\nu$ is self-adjoint,
\item[(iii)] $\Hi_\nu$ is an $L^2$ isometry, and
\item[(iv)] $\Hi_\nu A_\nu = \Omega^2 \Hi_\nu$.
\end{itemize}
The first of these is an immediate consequence of properties of the Fourier-Bessel integral defining $\Hi_\nu$,
but a proof may be found in~\cite{PST1} along with a proof of the fourth.
The second is obvious from the definition, and the third follows from the
first and second.

We define fractional powers of $A_\nu$ using the fourth property above:
\begin{displaymath}
A^{\sigma/2}_\nu = \Hi_\nu \Omega^{\sigma} \Hi_\nu.
\end{displaymath}
An integral kernel for $A^{\sigma/2}_\nu$ is given in~\cite{PST1},
\begin{equation}\label{intker}
(A^{\sigma/2}_\nu \phi)(r,\theta)
= \int^\infty_0 k^{\sigma}_{\nu,\nu}(r,s)\phi(s,\theta) s^{n-1}\,ds
\end{equation}
where
\begin{displaymath}
k^\sigma_{\nu,\nu} (r,s) = \begin{cases}{
\frac{2^{\sigma+1}\Gamma(\nu+\frac{\sigma}{2}+1)}{\Gamma(-\frac{\sigma}{2})\Gamma(\nu+1)}
\frac{s^{\nu-\lambda}}{r^{\sigma+\lambda+\nu+2}}
F(\nu+\frac{\sigma}{2}+1,\frac{\sigma}{2}+1;\nu+1;\frac{s^2}{r^2})}
& \mbox{ \,if\,\,\,} s<r,\cr
\frac{2^{\sigma+1}
\Gamma(\nu+\frac{\sigma}{2}+1)}{\Gamma(-\frac{\sigma}{2})\Gamma(\nu+1)}
\frac{r^{\nu-\lambda}}{s^{\sigma+\lambda+\nu+2}}
F(\nu+\frac{\sigma}{2}+1,\frac{\sigma}{2}+1;\nu+1;\frac{r^2}{s^2})
&  \mbox{ \,if\,\,\,}s>r.\end{cases}
\end{displaymath}
Here $F$ is the hypergeometric function 
\begin{equation}\label{def:F}
F(\alpha,\beta;\gamma;z) =
1+\frac{\alpha\cdot\beta}{1\cdot\gamma}z+\frac{\alpha(\alpha+1)\beta(\beta+1)}{1\cdot 2\cdot \gamma (\gamma+1)}z^2 + \dots
\end{equation}
  The integral in (\ref{intker}) may be interpreted in the usual Lebesgue sense for $\sigma < 0$.
See \cite{PST1} for the correct interpretation when $\sigma \ge 0$.

Hankel transforming both sides of the estimate we are trying to prove, we
see that we are reduced to proving
\begin{displaymath}
\|A_\nu^{-1/4-\alpha}\Omega^{1/2-2\alpha}\Hi_\nu S_\nu f\|_{L^2(\RR^{n+1})} \le C
\|\Hi_\nu f\|_{L^2(\RR^n)}
\end{displaymath}
where $\Hi_\nu S_\nu f$ solves
\begin{displaymath}
i\partial_t \Hi_\nu S_\nu f - \Omega^2 \Hi_\nu S_\nu f = 0, \qquad (\Hi_\nu S_\nu f)(0,\xi)
= (\Hi_\nu f)(\xi).
\end{displaymath}
But the solution to this initial value problem is just
\begin{displaymath}
(\Hi_\nu S_\nu f)(t,\xi) = \exp(-it|\xi|^2) (\Hi_\nu f)(\xi).
\end{displaymath}
Let $\Ft$ denote 
the Fourier transform in the $t$ variable,
\[
\Ft f(\tau,x) = \frac{1}{\sqrt{2\pi}}\int e^{-it\tau} f(t,x) dt.
\]
It is
an isometry of $L^2(\RR^{n+1})$ and it commutes with both
$A_\nu^{-1/4-\alpha}$ and $\Omega^{1/2-2\alpha}$, since these are defined
pointwise in~$t$.  Thus
\begin{displaymath}
\|A_\nu^{-1/4-\alpha}\Omega^{1/2-2\alpha}\Hi_\nu S_\nu f\|_{L^2(\RR^{n+1})}
=
\|A_\nu^{-1/4-\alpha}\Omega^{1/2-2\alpha}\Ft{\Hi_\nu S_\nu f}\|_{L^2(\RR^{n+1})}.
\end{displaymath}
From the calculation of the last paragraph we see that
\begin{displaymath}
(\Ft{\Hi_\nu S_\nu f})(\tau,\xi) = (\Hi_\nu f)(\xi) \delta(\tau-|\xi|^2)
\end{displaymath}
and hence
\begin{eqnarray*}
\lefteqn{(A_\nu^{-1/4-\alpha}\Omega^{1/2-2\alpha}\Ft{\Hi_\nu S_\nu f})(\tau,\xi)}\hspace{.5in}&&\\
& =& \int k^{-1/2-2\alpha}_{\nu,\nu}(|\xi|,s)s^{1/2-2\alpha} \delta(\tau-s^2)\Hi_\nu f(s\xi/|\xi|) s^{n-1}ds \\
& = & \frac{1}{2}
\tau^{\lambda-\alpha+1/4}k^{-1/2-2\alpha}_{\nu,\nu}
(|\xi|,\sqrt{\tau})(\Hi_\nu f)(\sqrt{\tau}\xi/|\xi|).
\end{eqnarray*}
We now compute the square of the $L^2(\RR^{n+1})$ norm of this quantity.
We square the absolute value and integrate over $\RR^{n+1}$,
replacing the Cartesian
coordinates $\xi$ with spherical coordinates $\rho$, $\theta$ to obtain
\begin{displaymath}
\frac{1}{4}\int^\infty_0\int^\infty_0\int_{S^{n-1}}
\tau^{2\lambda-2\alpha+1/2} (k^{-1/2-2\alpha}_{\nu,\nu}
(\rho,\sqrt{\tau}))^2|(\Hi_\nu f)(\sqrt{\tau}\theta)|^2
\,d\theta\,\rho^{n-1}d\rho\,d\tau
\end{displaymath}
or, making the change of variable $\omega=\sqrt{\tau}$,
\begin{displaymath}
\frac{1}{2}\int^\infty_0\int^\infty_0\int_{S^{n-1}}
\omega^{4\lambda-4\alpha+2} (k^{-1/2-2\alpha}_{\nu,\nu}
(\rho,\omega))^2|(\Hi_\nu f)(\omega\theta)|^2
\,d\theta\,\rho^{n-1}d\rho\,d\omega.
\end{displaymath}
Since $A^{-1/2-2\alpha}_\nu=A^{-1/4-\alpha}_\nu A^{-1/4-\alpha}_\nu$ we have
\begin{displaymath}
k^{-1-4\alpha}_{\nu,\nu}(r,t) = \int^\infty_0 k^{-1/2-2\alpha}_{\nu,\nu}
(r,s) k^{-1/2-2\alpha}_{\nu,\nu} (s,t)\,s^{n-1}ds.
\end{displaymath}
We apply this with $r=t=\omega$ and $s=\rho$ to evaluate the integral over
$\rho$, since $k(r,s)=k(s,r)$, obtaining,
\begin{displaymath}
\frac{1}{2}\int^\infty_0\int_{S^{n-1}}
\omega^{4\lambda-4\alpha+2} k^{-1-4\alpha}_{\nu,\nu}(\omega,\omega)|(\Hi_\nu f)(\omega\theta)|^2
\,d\theta\,d\omega
\end{displaymath}
or, using the explicit formula for~$k_{\nu,\nu}^\sigma$ given above and Gauss's formula for the
value of the hypergeometric function at~$z=1$,
\begin{displaymath}
2^{1-4\alpha}\pi\frac{\Gamma(\nu-2\alpha+\frac{1}{2})\Gamma(4\alpha)}{
\Gamma(\nu+2\alpha+\frac{1}{2})\Gamma(2\alpha+\frac{1}{2})^2}
\int^\infty_0\int_{S^{n-1}}
\omega^{n-1} |(\Hi_\nu f)(\omega\theta)|^2d\theta\,d\omega.
\end{displaymath}
The double integral is just the square of the $L^2(\RR^n)$ norm of $\Hi_\nu f$ and therefore of
$f$.  Thus we see that
\begin{displaymath}
\|\Omega^{-1/2-2\alpha}A^{1/4-\alpha}_\nu S_\nu f\|_{L^2(\RR^{n+1})}
= C_{\nu,\alpha}
\|f\|_{L^2(\RR^n)}
\end{displaymath}
where
\begin{equation}\label{def:Cnual}
C_{\nu,\alpha} = 2^{1/2-2\alpha}\sqrt{\pi\frac{\Gamma(\nu-2\alpha+\frac{1}{2})\Gamma(4\alpha)}{
\Gamma(\nu+2\alpha+\frac{1}{2})\Gamma(2\alpha+\frac{1}{2})^2}}.
\end{equation}
As promised, one can see immediately from this formula that $C_{\nu,\alpha}$ is finite if $0<\alpha<\frac{1}{4}+\frac{\nu}{2}$, and 
that it
is a decreasing function of~$\nu$ and hence of~$l$. Thus to finish the proof of the Theorem, we expand $f$ in
spherical harmonics, $f = \sum_{l=d}^\infty f_l$, use the triangle inequality on the left and the $L^2$-orthogonality of 
spherical harmonic subspaces on the right, to obtain (\ref{eq:smoothS}), with the constant 
$C = \min_{l\geq d} C_{\nu_l,\alpha} = C_{\nu_d,\alpha}$.

\cqfd

\subsection{Generalized Morawetz estimate for the wave equation}

In this section we will obtain a weighted-$L^2$ estimate for the wave equation (\ref{eq:lwe}) which is analogous to the one
obtained above for the Schr\"odinger equation.
 The  estimate for the free
case $a=0$
 was proved in \cite{Hos97}.  The Morawetz estimate, which is the 
estimate for the free wave equation obtainable from the Morawetz radial identity \cite{Mor68}
is a special case of our estimate below (see Corollary~\ref{twabis}.)
 We note that the simple proof of (\ref{eq:smoothS}) from \cite{BAK} applies as well to the free wave equation.
\begin{theorem}
\label{twa}
   Let $n\geq 2$ and $d\geq d_0(n)$ be integers,  let $0<\alpha<\frac{1}{4} + \frac{1}{2}\nu_d$
 and let $u$ be the solution
 to {\rm(\ref{eq:lwe})}.  There exists  a constant $C>0$ depending on $n,a,d,\alpha$ such that for all 
$f \in \dot{H}^{1/2}_{\geq d}(\RR^n)$
and $g \in \dot{H}^{-1/2}_{\geq d}(\RR^n)$,
   \begin{equation}
     \label{eq:ff666}
      \| \Omega^{-1/2-2\alpha} P_a^{1/4-\alpha} u\|_{L^2(\RR^{n+1})} \leq
 C( \|f\|_{\dot{H}^{\frac{1}{2}}(\RR^n)}+\|g\|_{\dot{H}^{-\frac{1}{2}}(\RR^n)}).
   \end{equation}
 \end{theorem}
{\em Proof}: Once again, we can work one spherical harmonic at a time, and thus we are solving
\[
\partial_t^2 u_l + A_\nu u_l = 0,\qquad u_l(0,x) = f_l(x),\ \partial_tu_l(0,x) = g_l(x)
\]
with $\nu>0$ as before.
Applying the Hankel transform (and suppressing the $l$ subscripts) we obtain the solution as
\[
\Hi_\nu u (t,\xi) = \cos( t|\xi|)\Hi_\nu f(\xi) + \frac{\sin(t|\xi|)}{|\xi|} \Hi_\nu g(\xi),
\]
and  Fourier-transforming in time we have
\[
\Ft \Hi_\nu u(\tau,\xi) = \frac{1}{\sqrt{|\xi|}}(\delta(\tau+|\xi|)h_+(\xi) +\delta(\tau-|\xi|)h_-(\xi))
\]
where 
\[
h_\pm(\xi) = \frac{1}{2} (\sqrt{|\xi|}\Hi_\nu f(\xi) \pm \frac{1}{i\sqrt{|\xi|}} \Hi_\nu g(\xi))
\]
so that
\begin{equation}\label{eq:norms}
\|h_\pm\|_{L^2(\RR^n)} \leq C (\|A_\nu^{1/4} f\|_{L^2(\RR^n)} + \|A_\nu^{-1/4}g\|_{L^2(\RR^n)}).
\end{equation}
Thus for $\tau>0$,
\[
(A_\nu^{-1/4-\alpha}\Omega^{1/2-2\alpha} \Ft \Hi_\nu u)(\tau,\xi) = \tau^{n-1-2\alpha}k_{\nu,\nu}^{-1/2-2\alpha}(|\xi|,\tau)h_-(\tau \xi/|\xi|)
\]
and for $\tau<0$ we change $\tau$ to $-\tau$  and $h_-$ to $h_+$ in the right hand side above. By the same  calculation as in 
the Schr\"odinger case performed in the previous section we then have
\[
\|A_\nu^{-1/4-\alpha}\Omega^{1/2-2\alpha} \Ft \Hi_\nu u\|_{L^2(\RR^{n+1})} = C_{\nu,\alpha}\|h_\pm\|_{L^2(\RR^n)}
\]
with $C_{\nu,\alpha}$ as in (\ref{def:Cnual}).  

In light of (\ref{eq:norms}) we have the desired estimate provided we can show that the Sobolev norms based on 
(small) powers of the operator $P_a$ are equivalent to standard Sobolev norms based on the Laplacian.  This can 
be accomplished in several ways. Here we will give a simple proof of a such a result based on Hardy's inequality.  (The 
range of $s$ here is not optimal, one could improve on it by using the more sophisticated machinery developed in Section \ref{sec:derivatives}).
\begin{prop}
\label{equiv-norm}
Let $n\geq 3$, $a+\lambda^2>0$ and $-1\leq s\leq 1$.  There exists  constants $C_1,C_2>0$ depending on $n,a,s$ such that
\begin{equation}\label{equiv}
C_1\|f\|_{\dot{H}^s(\RR^n)} \leq \|P_a^{s/2} f\|_{L^2(\RR^n)} \leq C_2 \|f\|_{\dot{H}^s(\RR^n)}
\end{equation}
for all $f \in \dot{H}^s(\RR^n)$.  For $n=2$ the same result holds for functions $f\in \dot{H}^s_{\geq 1}$.
\end{prop}

{\em Proof}: We use the following version of Hardy's inequality: For $n\geq 3$,
\begin{equation}\label{hardy}
\|\Omega^{-1}f\|_{L^2(\RR^n)}\| \leq \frac{1}{\lambda} \|f\|_{\dot{H}^1(\RR^n)}
\end{equation}
A similar inequality is easy to obtain for $n=2$ and $f$ orthogonal to radial functions: Let $f(r,\theta) = \sum_{n\ne 0} f_n(r) e^{in\theta}$
be the Fourier series expansion of such an $f$.  Then
\[
\int_0^{2\pi} |f(r,\theta)|^2 d\theta
 =  \sum_{n\ne 0} |f_n|^2 \leq \sum_{n\ne 0} n^2 |f_n|^2 = \int_0^{2\pi} |\partial_\theta f|^2 d\theta
\]
and thus
\[
\|\Omega^{-1}f\|_{L^2}^2 \leq \int_0^\infty \int_0^{2\pi} \frac{1}{r^2}|\partial_\theta f(r,\theta)|^2 d\theta\ rdr
\leq \|\nabla f\|_{L^2}^2
\]
which establishes the Hardy inequality in 2 dimensions, with constant = 1.  
Denoting now the $L^2(\RR^n)$ inner product by $\langle,\rangle$, we have
\[
\|P_a^{1/2} f\|_{L^2}^2 = \langle f, -\Delta f \rangle + a \|\Omega^{-1}f\|_{L^2}^2 
\]
Thus for $a>0$ by (\ref{hardy}),
\[
\| \nabla f\|_{L^2}^2 \leq \|P_a^{1/2} f\|_{L^2}^2 \leq (1+\frac{a}{\lambda^2}) \|\nabla f\|_{L^2}^2
\]
while for $a<0$ we have the same as above with the two inequality signs reversed.  Thus 
we obtain (\ref{equiv}) for $n\geq 3$ and $s=1$ with
$C_1 = \min\{\nu_0/\lambda,1\}$ and $C_2 = \max\{\nu_0/\lambda,1\}$.  For $n=2$ we have $C_1 =1$ 
and $C_2 = 1+a$. By duality, we obtain (\ref{equiv}) for $s=-1$ as well, and interpolating
between the  two endpoints establishes the claim.\cqfd
\section{Strichartz estimates}
\subsection{The Schr\"odinger equation}
\label{StriSchrod}
Here we prove that the solution to (\ref{eq:lse}) satisfies the same set of estimates as
that of the free Schr\"odinger equation, for $n\geq 2$:

\begin{theorem}\label{thm:strich}
  Let $f \in L^2$, $p\geq 2,q$ such that
  \begin{equation}\label{schrodadmis}
\frac{2}{p}+\frac{n}{q}=\frac{n}{2},\qquad (n,p)\ne(2,2).
\end{equation}
Let $u$ be the unique solution of (\ref{eq:lse}).  Then, provided  $a+\lambda^2>0$,
there exists a constant $C>0$ depending on $n,p,a$, such that
  \begin{equation}
    \label{eq:strichartz}
    \|u\|_{L^p_t(L^q_x)} \leq C \|f\|_{L^2}.
  \end{equation}
\end{theorem}
We will follow the strategy from \cite{RodS} (notice however that using the
end-point allows to shorten the argument for $n\geq 3$, as well as to
recover said end-point). We consider the
potential term as a source term, 
\begin{equation}
  \label{eq:ff1}
  i\partial_t u+\Delta u=\frac{a}{|x|^2} u,
\end{equation}
and integrate using $S(t)=\exp(it\Delta)$, the free evolution, to get
\begin{equation}
  \label{eq:duhamel}
  u(t)=S(t)f-i a\int_0^t S(t-s) \Omega^{-2} u(s) ds.
\end{equation}
The first term can be ignored, and we focus on the Duhamel term. We
postpone the $n=2$ case, and assume $n\geq 3$. Then, for the free
evolution, one has Strichartz estimates up to the end-point, namely the
pair $(p,q)=(2,\frac{2n}{n-2})$. We recall that these Strichartz
estimates hold in a slightly relaxed setting (\cite{KeeTao98}), 
\begin{equation}
  \label{eq:yeah}
\|  \int_0^t S(t-s) F(x,s)ds \|_{L^2_t(L_x^{\frac{2n}{n-2},2})} \leq C \| F \|_{L^2_t(L_x^{\frac{2n}{n+2},2})} ,
\end{equation}
where $L^{\alpha,\beta}$ are Lorentz spaces. We also note that the estimate (\ref{eq:strichartz})
for general $(p,q)$ satisfying (\ref{schrodadmis}) follows by interpolating between
 the above endpoint estimate and the well-known estimate $p=\infty,q=2$ corresponding to the conservation of charge
for the Schr\"odinger equation.

 Hence to prove our
estimate, all we need to check is $F= \Omega^{-2} u \in
L^2_t(L_x^{\frac{2n}{n+2},2})$. We have from Theorem~\ref{thm:smoothing}
with $d=0$, $\alpha = 1/4$ that
\[
\| F\|_{L_t^2L_x^{\frac{2n}{n+2},2}} = \| \Omega^{-1}
\Omega^{-1} u \|_{L_t^2L_x^{\frac{2n}{n+2},2}}  \leq  \| \frac{1}{|x|}\|_{L_x^{n,\infty}} \|
\Omega^{-1}u\|_{L_t^2L_x^2}  \leq  C \|f\|_{L^2},
\]
where we have made use of the generalized H\"older inequality \cite{One63}.  
This ends the proof for large dimensions.

In the $2D$ case, one needs to follow \cite{RodS} more closely, and
resort to the following lemma, proved in \cite{CK}. 
\begin{lemma}\label{lem:CK}
 Let $X,Y$ be two Banach spaces and let $T$ be a bounded linear operator from $L^\beta(\RR^+;X)$ to
 $L^\gamma(\RR^+;Y)$, such that $Tf(t)=\int_0^\infty K(t,s)f(s)ds$. Then the
 operator $\tilde{T} f(t)=\int_0^t K(t,s)f(s)ds$ is bounded from
 $L^\beta(\RR^+;X)$ to $L^\gamma(\RR^+;Y)$ when $\beta<\gamma$, and $\|\tilde{T}\| \leq c_{\beta,\gamma} \|T\|$ with
 $c_{\beta,\gamma} = (1-2^{1/\gamma-1/\beta})^{-1}$. 
\end{lemma}
 Using this lemma, one may forget
about the $\int_0^t$ in the Duhamel formula (\ref{eq:duhamel}), and replace it with an integral over all
times. Thus we define, for $h \in L^2(\RR^{2+1})$,
$$
Th(t)=S(t)\int S(-s) \Omega^{-1}h(s) ds.
$$
Consider now the estimate dual to (\ref{eq:smoothS}), with $d=1$, $a=0$ and
$\alpha = 1/4$.  It reads
\[
\| \int S(s) F(s) ds \|_{L^2_x} \leq C \|\ \Omega F \|_{L^2_{t,x}}
\]
Combining this with the Strichartz estimate for the free Schr\"odinger group $S(t)$, we see that 
$T: L^2(\RR,L^2) \to L^p(\RR,L^q)$ with $p,q$ as in the statement of the
Theorem. Since $p>2$, by the above Lemma~\ref{lem:CK}, the operator $\tilde{T}$ is
bounded on the same spaces.  On the other hand, from (\ref{eq:duhamel})
we have that, for $u$ the solution to (\ref{eq:lse}),
\[
u(t) = S(t) f - i a \tilde{T} (\Omega^{-1} u).
\]
We now use (\ref{eq:smoothS}) again with $\alpha = 1/4$, to conclude
\[
\|u\|_{L^p_tL^q_x} \leq C (1 + |a| C) \|f\|_{L^2}.
\]
Since we are in two space dimensions, we need to assume that the data
$f$ is orthogonal to radial functions, in order for (\ref{eq:smoothS})
to hold in the case $a=0,\alpha=1/4$.  On the other hand, Strichartz estimates for  (\ref{eq:lse})
 in the case of radial data are obtainable from 
the  case $a=0$ using the
conjugation procedure presented in \cite{PST1}. 
For future reference, we state  a more general result here:  
\begin{theorem}\label{thm:strichschrodrad} Let $f\in \dot{H}^s_{<d}(\RR^n)$.  Then there is a 
constant $C>0$ (depending on $n,d,p,s$)
such that $u$, the solution to
(\ref{eq:lse}) satisfies 
\begin{equation}\label{radialstrich}
\|(-\Delta)^{s/2}u\|_{L^p_tL^q_x} \leq C \|f\|_{\dot{H}^s}
\end{equation}
for all $(p,q)$ as in {\rm(\ref{schrodadmis})} and all $s$ such that
\[
-\min\{1+\sqrt{a+\lambda^2},\frac{n}{2} + \frac{2}{p}\} < s < \min\{  1+\sqrt{a+\lambda^2}-\frac{2}{p}, \frac{n}{2} \}
\]
\end{theorem}

{\em Proof:} Let $f = \sum_{l=0}^{d-1} f_l$ be the spherical harmonics decomposition of $f$, and let $u= \sum_{l=0}^{d-1}$ be the 
corresponding decomposition of the solution $u$.  On the $l$'th harmonic subspace, $-\Delta = A_\mu$ and 
$-\Delta + a|x|^{-2} = A_\nu$,
with $\mu = \lambda+l$, $ \nu = \nu_l(n,a)$.  We thus have 
\[
i\partial_t u_l - A_\nu u_l = 0, \qquad u_l(0) = f_l,\qquad l=0,\dots,d-1
\]
In \cite{PST1} it was shown that the operator $\K_{\mu,\nu}^0:= \Hi_\mu \Hi_\nu$ is a conjugation operator between $A_\mu$ and $A_\nu$, i.e. 
$A_\mu \K^0_{\mu,\nu} = \K^0_{\mu,\nu} A_\nu$.  Thus, to obtain the estimate (\ref{radialstrich}) for each $u_l$ by conjugation 
from the corresponding estimate in the case $a=0$, all we need to know 
are the continuity properties of $\K_{\mu,\nu}^0$ on appropriate Sobolev spaces.  In particular, for the 
above estimate (\ref{radialstrich}) we need
continuity of $\K^0_{\mu,\nu}$ on $\dot{H}^s$ for the right-hand side and continuity of $\K^0_{\nu,\mu}$ 
on $\dot{H}^s_q$ for the left.
It was shown in \cite{PST1} that $\K^0_{\mu,\nu}$ is continuous on $\dot{H}^s_r$ provided that
\[
-\min\{\lambda,\mu,\nu,\mu-s\} < \frac{n}{r}-\lambda < 2 +\min\{\lambda,\mu,\nu,\nu+s\}
\]
Applying this in  the two cases we need here gives the restriction on $s$ in the statement of the theorem.  We thus obtain the 
desired Strichartz estimate on each spherical harmonic subspace.  These need to be added up, which can be done since there
is only a finite number of them, and finally $\|f_l\|_{\dot{H}^s} \leq C_{l,s} \|f\|_{\dot{H}^s}$ where $C_{l,s}$ is the norm of 
the projection operator onto the $l$'th spherical harmonic subspace.
\cqfd

To conclude the proof of Theorem~\ref{thm:strich}, we use Theorem~\ref{thm:strichschrodrad}
 with $n=2$, $s=0$ and $d=1$ to obtain the estimate for the radial part of the data.  \cqfd

 \begin{remark}
\rm
   We end this section by some comments on dispersive estimates and
   their relationship to Strichartz estimates. The main goal in
   \cite{JouSofSog91} was to obtain the $L^\infty-L^1$ dispersive estimate,
   and they deduced Strichartz estimates as a corollary, applying the
    usual duality argument from the free case. However, it required
   considerably more work to obtain such a dispersive estimate in the
   presence of a potential, and it imposed the assumption of rapid decay on this
   potential. The arguments from \cite{RodS} bypass the dispersive
   estimate to prove Strichartz directly. In our setting, this
   type of approach is required since the $L^\infty-L^1$ dispersive estimate
   is known to fail (\cite{PST2}), at least for negative $a$. However,
   to perform the duality argument referred to in the above, one only needs an
   $L^{p'}-L^p$ dispersive estimate, for $p\sim 2n/(n-2)$. Hence,
   one may wonder whether such a restricted dispersive estimate holds true. 
   Estimates of this type have been obtained in \cite{GeoVis} for the wave equation
   with a potential decaying strictly faster than $|x|^{-2}$.  In the
   remaining part of this section we will indicate a simple way to get a
   dispersive estimate with a epsilon loss. The same argument would apply to the wave equation.
 \end{remark}
 We will prove the following result
\begin{prop}\label{prop:cute}
Let $n\geq 5$. Let $u$ be the unique solution of (\ref{eq:lse}). There exists a
  constant $C_\varepsilon$ (depending also on $n,a$), such that
  \begin{equation}
    \label{eq:dispSchrod}
    \|\Delta_j u(t,\cdot)\|_{L^{\frac{2n}{n-2}}} \leq C_\varepsilon
    \frac{2^{2j\varepsilon}}{t^{1-\varepsilon}} \|\Delta_j f\|_{L^{\frac{2n}{n+2}}},
  \end{equation}
where $\Delta_j$ is the usual frequency localization at $|\xi|\sim 2^j$.
\end{prop}
 Let us assume momentarily that we have obtained a
weighted $L^2$ dispersive estimate of Kato-Jensen type (\cite{KJ79}),
\begin{equation}
  \label{eq:katojensen}
  \int |u|^2\frac{dx}{|x|^{2\alpha}} \lesssim \frac{1}{t^{2\alpha}} \int
|f|^2{|x|^{2\alpha}}{dx},
\end{equation}
for some range of $\alpha$ (which will depend on the dimension). This estimate will play the role assigned to
(\ref{eq:smoothS}) in our proof of Strichartz.

We proceed with a nice observation we learned from J. Ginibre (\cite{Ginibre}).
 Recall
(\ref{eq:duhamel}),
$$
 u(t)=S(t)f-ia\int_0^t S(t-s) \Omega^{-2} u(s) ds,
$$
and write also the reversed Duhamel formula, where
$S_a(t)=\exp(-itP_a)$ plays the role of the free group.
Replace $u(s)$ in the first Duhamel formula by its expression coming
from the second, we get
\begin{eqnarray}
\label{dubduhamel}
  u(t) & = & S(t)f-i a\int_0^t S(t-s) \Omega^{-2} S(s) f ds \\
& & {}+a^2\int_0^t \int_0^s S(t-s) \Omega^{-2} S_a(s-\tau)
\Omega^{-2} S(\tau) f d\tau ds.\nonumber
\end{eqnarray}
Recall the dispersive estimate (with Lorentz spaces) for the free
group, which we state for a frequency localized data:
$1<p'\leq 2$,
$$
\|\Delta_j S(t)f\|_{L^{p,2}}\leq C\min(2^{nj(\frac{1}{p'}-\frac{1}{p})},
t^{-n(\frac{1}{2}-\frac{1}{p})})\|\Delta_j f\|_{L^{p',2}}.
$$
For the second term in (\ref{dubduhamel}) we use $p=2n/(n-2)$, and
compute explicitly the integral after using the free dispersion, to get  $\log(1+2^{2j}t)t^{-1}$. So we focus on the
third term. Freezing the time variables, we have, with
$p=\frac{2n}{n-2-2\mu}$, ($\mu$ should be thought as small)
\begin{eqnarray*}
 \| \Delta_j S(t-s) \Omega^{-2} S_a(s-\tau)
\Omega^{-2} S(\tau) \Delta_j f\|_{L^{p,2}}  \leq & \\
\min(2^{2j(1+\mu)},\frac{1}{(t-s)^{1+\mu}})
\frac{1}{(s-\tau)^{1-\mu}}\min(2^{2j(1+\mu)},\frac{1}{\tau^{1+\mu}})&\| 
 \Delta_j f\|_{L^{p',2}},
\end{eqnarray*}
where we have successively used free dispersion, generalized H\"older,
weighted $L^2$ dispersion (\ref{eq:katojensen}) with $\alpha=1-\mu$, generalized H\"older again, and finally
free dispersion. Evaluating the double time integral yields the
desired decay $2^{4\mu j}t^{-1+\mu}$, and interpolation with the $L^2$
bound gives the result (said $L^2$ bound resulting from
Strichartz which hold for both $S(t)$ and $S_a(t)$). 

We now prove (\ref{eq:katojensen}). In \cite{PST2} we obtained
some dispersive estimates in the radial case, but the argument would
equally apply to any spherical harmonic with uniform constants. Thus,
one could deduce (\ref{eq:katojensen}) from these estimates, using
generalized H\"older in the radial variable. However, we choose to
give a simple proof here, based on the following observation, which was
pointed out to us by I. Rodnianski: introduce $L=x/2+it\nabla$ the pseudo-conformal
vectorfield. $L$ commutes with the free evolution, but not with the
potential term. Introduce $C=L^2+at^2\Omega^{-2}$. Then a simple sequence of
computations gives
$$
[i\partial_t -P_a,C] u = it(2a\Omega^{-2}+a x\cdot \nabla \Omega^{-2}) u.
$$
Due to the special form of the potential, $2\Omega^{-2}+ x\cdot \nabla
\Omega^{-2}=0$, and $C$ commutes with the equation. Therefore,
$\|Cu\|^2_2$ is conserved. Now, another simple computation gives
 $\Delta M f=\frac{1}{t^2} M L^2 f$, where $Mf=e^{i\frac{x^2}{4s}} f$. Adding the potential,
$
P_a M f=-\frac{1}{t^2} M C f,
$
since obviously $M$ and $\Omega^{-2}$ commute.
Provided $a>-(n-2)^2/4$, we can apply Hardy to the
left-hand side (this will require $n\geq 5$), set $f=u$, to get
$$
\|\frac{Mu}{|x|^2}\|_2 \leq \frac{C}{t^2}\|M Cu\|_2,
$$
and the $M$ is now irrelevant, so that
$$
\|\frac{u}{|x|^2}\|_2 \leq \frac{C}{t^2}\| Cu_0\|_2
=\frac{C}{t^2}\| |x|^2 u_0\|_2.
$$
Interpolating this with $L^2$ conservation for $u$ gives the desired
estimate. \cqfd

\begin{remark}\rm
  We note that such a weighted $L^2$ dispersive estimate (with a decay
  greater than $t^{-1}$) is all which is required to run the
  perturbative argument from \cite{RodS}, with the free evolution
  replaced by $S_a(t)$. In fact, all the other estimates which are
  required are contained in Theorem \ref{thm:smoothing}. Hence, one
  could in principle obtain Strichartz estimates for (non radial) perturbations of
  the inverse square potential which decay $\varepsilon$ faster.
\end{remark}

\subsection{The wave equation}
In this section we will obtain Strichartz estimates for the operator
$\Box_a=\partial^2_t+P_a$ from the generalized Morawetz estimate~(\ref{eq:ff666})
 as we did for the Schr\"odinger
equation. We  present two different results: the first one
is in some sense the true equivalent of Theorem \ref{thm:smoothing}; it
uses Lemma \ref{lem:CK} to bypass issues related to the varying
degrees of smoothness in Strichartz estimates for the wave
equation. The second result shows one can in fact recover the
end-point Strichartz estimate if needed. These two 
results will merge into a  theorem on generalized Strichartz estimates with derivatives, in
the last section.
\begin{theorem}\label{thm:strichwave}
  Let $u$ be the solution to {\rm(\ref{eq:lwe})} with Cauchy data $(f,g)
  \in \dot{H}^{\frac{1}{2}}\times \dot{H}^{-\frac{1}{2}}$.  Let  $p\geq 2$, and $q$ be
  such that
  $\frac{2}{p}+\frac{n-1}{q}=\frac{n-1}{2}$ ($p>2$ if $n=3$ and $p>4$ if $n=2$). Then,
  provided  $a+\lambda^2>0$,
  \begin{equation}
    \label{eq:strichartzwave}
    \| (-\Delta)^{\sigma/2} u\|_{L^p_t(L^q_x)} \leq 
C (\|f\|_{\dot{H}^{\frac{1}{2}}}+\|g\|_{\dot{H}^{-\frac{1}{2}}}).
  \end{equation}
where $\sigma=\frac{1}{p}+\frac{n}{q}-\frac{n-1}{2}$ (gap condition).
\end{theorem}
One could of course replace the norm on the left with the appropriate
Besov norm.

{\em Proof}: We write the solution to the wave equation with a potential (\ref{eq:lwe}) as the sum
of the linear solution to the free wave equation 
$$
\dot{W}(t)f+W(t)g=\mathcal{F}^{-1}(\cos(t|\xi|)\mathcal{F}{f}+\frac{\sin(t|\xi|)}{|\xi|} \mathcal{F}{g}),
$$
where $\mathcal{F}$ is the Fourier transform in the space variable $x$,  plus a Duhamel term
\begin{equation}
  \label{eq:ff10}
  a\int_0^t W(t-s) \Omega^{-2} u(s)ds.
\end{equation}
Since $W(t-s) = -\dot{W}(t) {W}(s) + {W}(t)\dot{W}(s)$, we obtain two terms in the above.  We will deal with the 
first one, the treatment of the second term being similar. We are going to use Lemma~\ref{lem:CK}, thus we set
$$
T h(t) :=\dot{W}(t)\int W(s) \Omega^{-1} h(s) ds.
$$
Once again, using the Strichartz estimate for the free wave equation, combined with the dual to 
estimate (\ref{eq:ff666}) (for $n\geq 2$, $d=d_0(n)$, $a=0$ and $\alpha = 1/4$),
\begin{eqnarray*}
\|Th\|_{L^p\dot{H}^\sigma_q} & \leq & C \|\int W(s) \Omega^{-1} h(s) ds\|_{\dot{H}^{1/2}(\RR^n)}\\
& \leq &C \|h\|_{L^2(\RR^{n+1})}
\end{eqnarray*}
with $p$, $q$ and $\sigma$ as in the statement of the Theorem.  By Lemma~\ref{lem:CK}, the corresponding operator
$\tilde{T}$ satisfies the same estimate as $T$ (with a different constant).  On the other hand, the solution to
(\ref{eq:lwe}) is
\[
u(t) = \dot{W}(t)f+W(t)g + a \tilde{T}(\Omega^{-1} u)
\]
thus using (\ref{eq:ff666}) one more time we obtain the desired result.  For $n=2$ we need to assume that the 
data $f,g$ are orthogonal to radial functions, for (\ref{eq:ff666}) to hold in the $a=0,\alpha=1/4$ case.
On the other hand, Strichartz estimates for the wave equation (\ref{eq:lwe}) 
in case of radial data
 were proven in \cite{PST1} using the conjugation method. For future reference we
 quote a more general result here. We note that the restriction on $\gamma$ and hence on $\sigma$, which is the
number of derivatives that can be taken, comes from the requirement of the continuity of the 
conjugation operator $\K^0_{\mu,\nu}$ on the spaces involved (see \cite{PST1} for details).
\begin{theorem}\label{thm:strichwaverad}
Let $n\geq 2$,  $2\leq q <\infty$ and let $p,\gamma,\sigma$ be such that 
\begin{equation}\label{waveadmis}
\frac{1}{p} \leq \min\{\frac{1}{2},\frac{n-1}{2}(\frac{1}{2}-\frac{1}{q})\},
\qquad \sigma = \gamma+\frac{1}{p} - n(\frac{1}{2}-\frac{1}{q}).
\end{equation}
 For integer $d\geq 1$ let
$f\in \dot{H}^{\gamma}_{<d}$ and $g \in \dot{H}^{\gamma-1}_{<d}$.  Then there exists a constant $C>0$ depending
on $n,q,p,\gamma,d$ such that the solution to {\rm(\ref{eq:lwe})}
satisfies 
\begin{equation}\label{radstrichwav}
\|u\|_{L^p\dot{H}^\sigma_q} \leq C(\|f\|_{\dot{H}^{\gamma}}+\|g\|_{\dot{H}^{\gamma-1}}).
\end{equation}
provided
\[
-1-\nu_0<\gamma <1+\nu_0 - \frac{1}{p}
\]
\end{theorem}
To complete the proof of Theorem~\ref{thm:strichwave}, we use the above estimate with $n=2$, $d=1$, and $\gamma=1/2$
for the radial part of the data.
\cqfd

We now turn to an endpoint estimate:
\begin{theorem}\label{thm:strichwavebis}
  Let $n\geq 4$, let 
\begin{equation}\label{magicnumbers}
(p,q,\gamma,\sigma) = (2,\frac{2(n-1)}{n-3},\frac{n-3}{2(n-1)},\frac{-2}{n-1}),
\end{equation}
 and let
 $u$ be the solution to {\rm(\ref{eq:lwe})}. Then, provided  $a+\lambda^2>0$,
  \begin{equation}
    \label{eq:strichartzwavebis}
    \| (-\Delta)^{-\sigma/2} u\|_{L^{p}_t(L^q_x)} 
\leq C (\|f\|_{\dot{H}^\gamma}+\|g\|_{\dot{H}^{\gamma-1}}).
  \end{equation}
\end{theorem}
The  exponents here are chosen so as to postpone until Section~\ref{sec:derivatives} the
unpleasant issues related to commuting the free Laplacian with its
counterpart with potential.

{\em Proof of Theorem~\ref{thm:strichwavebis}}: We start with a simple corollary of Theorem \ref{twa}. 
Using the fact
that $P_a$ commutes with the solution,  and the equivalence
(\ref{equiv}), we have
\begin{corol}
\label{twabis}
   Let $n,d,\alpha,u$ be as in Theorem~\ref{twa}. Then
   \begin{equation}
     \label{eq:ff666bis}
      \| \Omega^{-1/2-2\alpha} u\|_{L^2(\RR^{n+1})} \leq
 C( \| f\|_{\dot{H}^{2\alpha}(\RR^n)}+\|g\|_{\dot{H}^{2\alpha-1}(\RR^n)}).
   \end{equation}
\end{corol}
We note that the case $a=0$, $\alpha=\frac{1}{2}$ is the Morawetz estimate \cite{Mor68}, which states that for 
the solution $u$ of the linear wave equation $\Box u = 0$ in four or more dimensions,
$\iint u^2/|x|^3 dxdt$ is bounded by the energy of the initial data.

Using the above estimate~(\ref{eq:ff666bis}), one can proceed as we did for Schr\"odinger on
the Duhamel term (\ref{eq:ff10}). Take $2\alpha=\frac{n-3}{2(n-1)}$,
then
$$
\Omega^{\frac{1}{n-1}-1} u \in L^2_t(L^2_x).
$$
From this, writing
$$
\Omega^{-2} u = {\Omega^{-(1+\frac{1}{n-1})}}
  {\Omega^{\frac{1}{n-1}-1}} u,
$$
and using $|x|^{-1-\frac{1}{n-1}} \in L^{n-1,\infty}_x$, the generalized
H\"older inequality for Lorentz spaces yields
$$
\Omega^{-2} u \in L^2_t(L^{2\frac{n-1}{n+1},2}_x).
$$
We can then finish the proof by applying the usual end-point to
end-point Strichartz estimate,
\begin{equation}\label{eq:es}
\| |\nabla|^{-\frac{2}{n-1}} \Box^{-1}
F\|_{L^2_t(L^{2\frac{n-1}{n-3},2}_x)} \leq C \|F\|_{L^2_t(L^{2\frac{n-1}{n+1},2}_x)}.
\end{equation}
\cqfd
\section{Estimates with derivatives}\label{sec:derivatives}
When dealing with nonlinear applications for either the Schr\"odinger or
the wave equation, the $f\in L^2$ or $\dot{H}^{1/2}$ estimates obtained above are often not
enough, and one needs to consider estimates involving derivatives of the
solution and the data. Moreover, the correct number of derivatives needed could be fractional.  
For the wave equation this is the case even at
the linear level, as we already saw when trying to prove Strichartz
estimates, where spaces such as $\dot{H}^{1/2}$ and $\dot{H}^{-1/2}$ naturally
appear. In the free case, frequency localizations commute with the flow,
hence one may immediately deduce estimates for $f\in \dot{H}^s$ from an
$L^2$ or $\dot{H}^{\frac{1}{2}}$ estimate, through Littlewood-Paley. In
our setting, this is no longer true however, and one needs to replace
frequency localizations $\phi(\sqrt{-\Delta})$ by those based on the
operator $P_a$. At the same time we need the final estimate to be
phrased in terms of standard Sobolev spaces, based on powers of
$-\Delta$, otherwise the estimate would be useless in nonlinear
applications (unless one studies carefully the multiplication properties
of spaces based on the operator $P_a$). At issue is therefore the lack of
commutation between the two localizations, which we are now going to
address.  Here we will obtain various estimates in weighted spaces   for the products
of the two projections, which we will then use to deduce generalized Strichartz estimates with derivatives
from the same weighted-$L^2$ estimates obtained above.  We will end this Section with a nonlinear application
that illustrates the need for derivative estimates.

Let $\Delta_j$ be the usual dyadic frequency localization at
  $|\xi|\sim 2^j$, and let $\Pi_k$ be the localization with respect to $\sqrt{P_a}$.  More precisely,
let $\mathbf{\beta}_0\in C^\infty_0(\RR^+)$ denote the standard bump function supported in  $[\frac{1}{2},2]$, with the property that
$\sum_j (\mathbf{\beta}_0(2^{-j}x))^2 = 1$ for all $x\in\RR^+$.  
and let $\beta_j(\xi) := \mathbf{\beta}_0(2^{-j}|\xi|)$.  Let $\Hi_\nu$ denote the Hankel transform of order $\nu$.  Let $\Delta_j^l$ and 
$\Pi_k^l$ denote the restrictions to $L^2_l$, the $l$-th spherical harmonic subspace of $L^2$, of the above projections $\Delta_j$ and 
$\Pi_k$ respectively.  It was shown in \cite{PST2} that
\[
\Delta_j^l = \Hi_\mu \beta_j \Hi_\mu,\qquad \Pi_k^l = \Hi_\nu\beta_k
\Hi_\nu
\]
where $\mu = \lambda+l$, $\lambda = \frac{n-2}{2}$, and $\nu = \sqrt{\mu^2+a}$.
Let us define the following operators:
\[
J_{jk} : = \Delta_j \Omega^{-2} \Delta_k, \qquad M_{kl} := \Delta_k \Pi_l,\qquad N_{lm} := \Pi_l \Delta_m.
\]
We are going to need
estimates on the operators $J_{jk},M_{kl},N_{lm}$, in weighted-$L^2$ spaces. These will be provided in the
 next three Lemmas.  The estimate needed for $N_{lm}$ is the 
almost orthogonality lemma for the projectors, the radial version of which was given in \cite{PST2}:
\begin{lemma}\label{lem:Mjk} Let $n\geq 2$, $d\geq 0$.  
For all positive $$\epsilon_1<\min\{\lambda+d,\nu_d\}+1,$$
 there exists a constants $C>0$ such that for all $j,k\in\ZZ$ and $f\in L^2_{\geq d}(\RR^n)$,
  \begin{equation}
    \label{eq:almost-orth}
    \|M_{jk} f\|_{L^2}, \| N_{jk} f\|_{L^2}\leq C 2^{-\epsilon_1|j-k|}\|f\|_{L^2}.
  \end{equation}
\end{lemma}
{\em Proof}:   Recalling the
definition of $\K_{\mu,\nu}^0$,
\[
\Delta_j^l\Pi_k^l = \Hi_\mu \beta_j \K_{\mu,\nu}^0 \beta_k \Hi_\nu.
\]
Since the Hankel transforms appearing at the extremes are $L^2$-isometries, the problem reduces to
showing that the operator 
\[
L_{jk} := \beta_j \K_{\mu,\nu}^0 \beta_k
\]
 is bounded on 
$L^2_l$, with a norm that is bounded independent of $l$. For $j$ close to $k$ this is obviously true by the 
boundedness of each factor. We note that here as well as
 in the Proposition that follows, only the support properties of $\beta_j$ are used,
and not that their squares form a partition of unity.   

Now for $j\ne k$,
\begin{equation}\label{Ljk}
L_{jk} f = \beta_j(r) \int_0^\infty k_{\mu,\nu}^0(r,s) \beta_k(s) f(s) s^{n-1} ds
\end{equation}
We  recall the formula for the integral kernel $k_{\mu,\nu}^0$
obtained in \cite{PST1}: 
\begin{equation}\label{kform}
k^0_{\mu,\nu}(r,s) =
\frac{2\Gamma(\frac{\mu+\nu}{2}+1)}{\Gamma(\frac{\nu-\mu}{2})\Gamma(\nu+1)}
\frac{s^{\nu-\lambda}}{r^{\lambda+\nu+2}}F(\frac{\mu+\nu}{2}+1,\frac{\nu-
\mu}{2}+1;\nu+1;(\frac{s}{r})^2)
\end{equation}
Here $F$ is the hypergeometric function defined by (\ref{def:F}).
 The
above formula for $k_{\mu,\nu}^0$ is valid for $s<r$.  For $s>r$ one
needs to switch $s$ and $r$, and switch $\mu$ and $\nu$ in the formula.
Turning now to $L_{jk}$ for $|j-k|\geq 3$, we see that in  (\ref{Ljk})
we have $s\sim 2^j$ and $r \sim 2^k$, thus
we either have $s\leq r/2$ or $s\geq 2r$.  Therefore the last argument of the
hypergeometric function in (\ref{kform}) will always be in $[0,\frac{1}{2}]$.
It is then easy to see from (\ref{def:F}) that $|F| \leq C$ independent of $l$.
Moreover, $\nu = \mu + O(\frac{1}{l})$ for large $l$, and from Stirling's formula, 
\[
\frac{2\Gamma(\frac{\mu+\nu}{2}+1)}{\Gamma(\frac{\nu-
\mu}{2})\Gamma(\nu+1)} = O(l^{-1+c/l})<C.
\]
By the procedure outlined in \cite[\S 3.1]{PST1}, the resulting pointwise bound for $k_{\mu,\nu}^0$  gives the
desired $L^2$ bound, namely,
\[
\|L_{jk} f\|_{L^2_0} \leq C  2^{-\delta|j-
k|}\|f\|_{L^2_0}.
\]
for any $\delta<\delta_l:=\min\{\mu,\nu\}+1$. Summing over $l$ and using orthogonality of spherical harmonics in $L^2$, 
we obtain (\ref{eq:almost-orth}), with $$\epsilon_1< \min_{l\geq d} \delta_l = \min\{\lambda+d,\nu_d\}+1.$$
\cqfd

Next we obtain a weighted-$L^2$ estimate for $M_{kl}$. 
 First we need the following general result:  
\begin{prop}\label{prop:neat} Let $n\geq 2$ and $d\geq d_0(n)$ be fixed integers.  
There exists a constant $C>0$ (depending on $n$ and $d$) such that for all
$f\in L^2_{\geq d}(\RR^n)$,
 \[
\|\beta_0 (-\Delta)^{-1} \beta_m f\|_{L^2} \leq C 2^{-|m|(\lambda+d)+m}\|f\|_{L^2}
\]
\end{prop}

{\em Proof}: Let us first assume $|m|\geq 3$. Let $K(x)$ denote the Newtonian potential, and let $T_d(x;y)$ be the Taylor 
polynomial of degree $d-1$ in $y$ for $K(x-y)$.  For example, for $n\geq 3$, we have that up to a constant factor,
\[
T_3(x;y) = |x|^{-n+2} + (n-2)|x|^{-n} x\cdot y + \frac{n-2}{2} |x|^{-n-2}\{n(x\cdot y)^2 -|x|^2 |y|^2\}.
\]
Let $K_d(x,y):= \max\{T_d(x;y),T_d(y;x)\}$ It is easy to see that $K_d \in L^2_{< d}$ and thus for $f\in L^2_{\geq d}$
we may write
\begin{eqnarray*}
(-\Delta)^{-1}\beta_m f(x)& = & \int K(x-y) \beta_m(y)f(y) dy\\
& = & \int (K(x-y) - K_d(x,y))\beta_m(y) f(y) dy.
\end{eqnarray*}
We have $|x|\sim 1$ and $|y|\sim 2^m$.  Since $|m|\geq 3$ we have that either
  $|y|\leq\frac{1}{2} |x|$ or $|y|\geq2|x|$. We can then 
check easily that
\begin{eqnarray*}
|K(x-y) - K_d(x,y)| &\leq & C
\frac{(\min\{|x|,|y|\})^d}{(\max\{|x|,|y|\})^{n+d-2}} \\
&\leq &C
\left\{\begin{array}{ll} 2^{-m(n+d-2)} & m\geq 3\\ 2^{md} & m\leq
-3\end{array}\right.
 \end{eqnarray*}
  We then have
\begin{eqnarray*}
\|\beta_0(x)(-\Delta)^{-1}\beta_m(y) f(y) \|_{L^2} & \leq & \sup|K(x-y) - K_d(x,y)|
 \|f\beta_m\|_{L^1} \\
& \leq & C \|f\|_{L^2} \left\{\begin{array}{ll}
2^{-m(\lambda+d-1)} & m\geq 3\\ 2^{m(\lambda+d+1)} & m\leq -3\end{array}
\right.
\end{eqnarray*}
which establishes the claim, for $|m|\geq 3$.  For $|m|\leq 2$ we simply observe that if we let
 $\tilde{K}(x,y) := \beta_0(x)K(x-y)\beta_m(y)$
then $\|\tilde{K}\|_{L^1L^\infty},\|\tilde{K}\|_{L^\infty L^1} < C$ and thus the operator corresponding to $\tilde{K}$ maps 
any $L^p$ into itself.\cqfd

\begin{lemma}
\label{x-Mx} Let $n\geq 2$, $d\geq d_0(n)+1$, and $0\leq\eta\leq 2$.  
There exist  a constant $C>0$ such that for all $j,k\in\ZZ$ and for 
all $f\in(L^2_{\geq d}(\RR^n)$,
\begin{equation}\label{est:Mjk1}
\|\ \Omega^{-\eta} M_{jk} \Omega^{\eta} \ f\|_{L^2} \leq C 2^{-\epsilon_2|j-k|}\|f\|_{L^2}
\end{equation}
for any $\epsilon_2<\min\{\lambda+d,\nu_d\}+1-\eta$.
\end{lemma}
{\em Proof}: We can obtain this estimate by interpolating between (\ref{eq:almost-orth}) and the following estimate
\begin{equation}\label{est:Mjk2}
\|\ \Omega^{-2} M_{jk} \Omega^2 f \|_{L^2} \leq C2^{-\gamma|j-k|} \|f\|_{L^2}
\end{equation}
for some appropriate $\gamma>0$.
To prove this estimate, or equivalently, its dual estimate, 
we note that on $L^2_l$,
\begin{eqnarray*}
\Omega^{2} \Pi_j \Delta_k \Omega^{-2} &=& \Omega^2 \Hi_\nu \beta_j \Hi_\nu \Hi_\mu \beta_k \Hi_\mu \Omega^{-2}\\
& = & \Hi_\nu A_\nu \beta_j \K_{\nu,\mu}^0 \beta_k A_\mu^{-1} \Hi_\mu,
\end{eqnarray*}
where the $A$'s are as in (\ref{def:Anu}).  Once again, the Hankel transforms at the two extremes being $L^2$-isometries, the above reduces to proving
an estimate for the operator $L_{jk}$ defined above, namely, 
\[
\|A_\nu L_{jk}A_\mu^{-1} f_l \|_{L^2} \leq C2^{-\gamma|j-k|}
 \|f_l\|_{L^2}
\]
We have $A_\nu = A_\mu + \frac{a}{r^2}$ and $A_\nu \K_{\nu,\mu}^0 = \K_{\nu,\mu}^0A_\mu$. Thus
\begin{eqnarray*} 
A_\nu L_{jk} A_\mu^{-1} f_l & = & (A_\mu \beta_j) \K_{\nu,\mu}^0 \beta_k A_\mu^{-1} f_l \\
&&\mbox{} + \beta_j \K_{\nu,\mu}^0 (A_\mu \beta_k) A_\mu^{-1} f_l \\
&& \mbox{} + \beta_j \K_{\nu,\mu}^0 \partial_s\beta_k \partial_s A_{\mu}^{-1} f_l \\
&& \mbox{} + \partial_r \beta_j (\partial_r \K_{\nu,\mu}^0) \beta_k A_\mu^{-1} f_l \\
&&\mbox{} + L_{jk} f_l \\
& := & I + II + III + IV + V.
\end{eqnarray*}
We have estimated $V$ already.  For the other four, first note that since what we want to prove, namely
(\ref{est:Mjk2}), is scale-invariant, in estimating any of the pieces $I$-$IV$ we can set either $j$ or $k$ equal to zero.
Thus for $I$ it is enough to show
\begin{equation}\label{est:I}
\sum_l\|(A_\mu\beta_j)\K_{\nu,\mu}^0\beta_0 A_\mu^{-1}f_l\|_{L^2} \leq C2^{-\gamma
j}\|f\|_{L^2}. \end{equation}
Now $A_\mu \beta_j(x) = 2^{-2j} \tilde{\beta}_j(x)$ where $\tilde{\beta}_j(x) = (A_\mu \beta_0)(2^{-j}x)$ is bounded independently of $j$ and
has the same support as $\beta_j$.
Using that $\sum_k \beta_k^2 = 1$ and that $\beta_j\beta_k = 0$ for $|j-k|>2$, we have
\begin{eqnarray*} 
\tilde{\beta}_j \K_{\nu,\mu}^0\beta_0 A_\mu^{-1} f_l & = &
 \tilde{\beta}_j \K_{\nu,\mu}^0\beta_0\sum_p \beta_p^2 A_\mu^{-1} \sum_m \beta_m^2 f_l\\
& = & \sum_{p=-1}^1 \sum_{m\in \ZZ} \tilde{\beta}_j \K_{\nu,\mu}^0\beta_0 \beta_p^2 A_\mu^{-1} \beta_m^2 f_l.
\end{eqnarray*}
It is enough to estimate the term with $p=0$, the other two being similar.  We then have
\[
\|\tilde{\beta}_j \K_{\nu,\mu}^0\beta_0 \beta_0^2 A_\mu^{-1} \beta_m^2 f_l\|_{L^2} \leq C 2^{-|j|\delta}\|\beta_0^2
 A_\mu^{-1} \beta_m^2 f_l\|_{L^2}
\]
by the previous Lemma, for any $\delta<\delta_l:=\min\{\mu,\nu\}+1$.  On the other hand, we have, 
by Proposition~\ref{prop:neat}, that
\begin{eqnarray*}
\sum_l \sum_m \|\beta_0^2 A_\mu^{-1} \beta_m^2 f_l\|_{L^2}& = & \sum_m \|\beta_0^2 \Delta^{-1} \beta_m^2 f \|_{L^2} \\
& \leq &  C \|f\|_{L^2} \sum_m 2^{-|m|(\lambda+d)+m} \\
& \leq & C \|f\|_{L^2}
\end{eqnarray*}
since by our assumptions $\lambda+d>1$, and thus
\[
\sum_{l=1}^\infty\|(A_\mu\beta_j)\K_{\nu,\mu}^0\beta_0 A_\mu^{-1}f_l\|_{L^2} \leq C
2^{-2j} 2^{-|j|\delta}\|f\|_{L^2} \leq C2^{-\gamma j}\|f\|_{L^2} \]
 for
\[
\gamma < \min_{l\geq d} \delta_l - 2 = \min\{\lambda+d,\nu_d\} -1.
\]
Estimating $II$ and $III$ is entirely analogous to the above, and for $IV$ we need to use the explicit form of the
kernel (\ref{kform}). Using the fact that $r$ and $s$ are well-separated, the series (\ref{def:F}) can be differentiated
term-by-term, thus obtaining corresponding decay rates for the derivative kernel, which in turn give the desired estimate
for $IV$ by the same procedure as above. 

We thus have the estimate (\ref{est:Mjk2}), with $\gamma$ as in the above.  Now,
 interpolating between (\ref{eq:almost-orth}) and (\ref{est:Mjk2}) gives the desired estimate 
(\ref{est:Mjk1}).
\cqfd

Finally, we need to estimate $J_{jk}$ on weighted-$L^2$ as well:
\begin{lemma}
\label{xJx}
 Let $n\geq 2$, $d\geq d_0(n)$ and $0\leq \zeta \leq 2$. There exist a constant $C>0$  such that for all
$f\in L^2_{\geq d}(\RR^n)$ and $j,k \in \ZZ$,
\begin{equation}\label{est:DjVDk}
\|\ \Omega^\zeta\ J_{jk}\ \Omega^{2-\zeta}\ f \|_{L^2} \leq C 2^{-
\epsilon_3 |j-k|}\|f\|_{L^2}
\end{equation}
for all $\epsilon_3 < \lambda+d-|1-\zeta|$.
\end{lemma}

{\em Proof}: After taking the Fourier transform, the estimate to prove is 
\[
\| \Delta_{\xi}^{\zeta/2} \beta_j \Delta_{\xi}^{-1} \beta_k \Delta_{\xi}^{1-\zeta/2} \hat{f} \|_{L^2} \leq C 2^{-
\epsilon_3|j-k|} \|\hat{f}\|_{L^2}
\]
where $-\Delta_{\xi} = \sum_i \partial^2/\partial \xi_i^2$  is the Laplacian in the Fourier variable.  We
will prove this by interpolation:  Let $T_{jk} := \beta_j \Delta_{\xi}^{-1} \beta_k$.
We then need to show that $T_{jk}$ maps $L^2$ into the homogeneous Sobolev space
$\dot{H}^2 $, and that it maps
$\dot{H}^{-2}$ into $L^2$.
Also note that $T_{jk} = T^*_{kj}$.  Thus it is enough to show that 
\begin{equation}\label{Tjk}
\|T_{jk} g\|_{\dot{H}^2} \leq C_{jk} \|g\|_{L^2}
\end{equation}
for some constant $C_{jk}$, in order to obtain via interpolation that 
\[
\|T_{jk} g\|_{\dot{H}^{\zeta}} \leq C_{jk}^{\zeta/2} C_{kj}^{1-\zeta/2} \|g\|_{\dot{H}^{\zeta-2}}
\]
The desired estimate would then follow by setting $g = \Delta_{\xi}^{1-\zeta/2}\hat{f}$. 

To prove (\ref{Tjk}) we again note that by scaling we can set $j=0$, and
estimate \[
\| \Delta_{\xi} \beta_0 \Delta_{\xi}^{-1} \beta_k g \|_{L^2}.
\]
We have
\begin{equation}\label{3terms}
\Delta_{\xi} \beta_0 \Delta_{\xi}^{-1} \beta_k g = \beta_0 \beta_k g + 2 \nabla \beta_0
\cdot \nabla \Delta_{\xi}^{-1} \beta_k g + (\Delta_{\xi}\beta_0 ) (\Delta_{\xi}^{-1}\beta_k g)
\end{equation}
The first term on the right is easy to estimate since $\beta_0\beta_k \equiv 0$ if
$|k|\geq 2$. For the third term, we use Proposition~\ref{prop:neat}, which
gives
\[
\|(\Delta_{\xi}\beta_0 ) (\Delta_{\xi}^{-1}\beta_k g)\|_{L^2} \leq C2^{-|k|(\lambda+d)+k}
\]

Similar argument applies to the middle term in (\ref{3terms}), and we thus obtain (\ref{Tjk}) with 
\[
C_{jk} := C2^{-|j-k|(\lambda+d)+j-k}
\]
Hence, 
$$
C_{jk}^{\zeta/2}C_{kj}^{1-\zeta/2} = 
C2^{-|j-k|(\lambda+d)-(j-k)(1-\zeta)} \leq C2^{-|j-k|(\lambda+d-|1-\zeta|)},
$$
 and we have the desired estimate for
$|j-k|\geq 3$.  For $j$ close to $k$, on the other hand, we can estimate the
two factors that make up $J_{jk}$ separately, i.e. it is enough to show that
$\Delta_{\xi}^{1/2} \beta_j\Delta_{\xi}^{-1/2}$ is bounded on $L^2$ independent of $j$.  By
interpolation, this further reduces to proving that $\Delta_{\xi} \beta_j \Delta_{\xi}^{-1}$ is
bounded on $L^2$, or equivalently, that multiplication by $\beta_j$ is bounded on
$\dot{H}^2$, which is easily seen to be the case.
\cqfd

\subsection{The Schr\"odinger equation}

To obtain a frequency-localized version of (\ref{eq:strichartz}), we recall that the solution
to (\ref{eq:lse}) can be expressed as
\[
u(t) = S(t) f -i a \int_0^t S(t-s) \Omega^{-2} S_a(s) f \ ds
\] 
where $S_a(s) = e^{-isP_a}$ is the  Schr\"odinger group associated to $P_a$.  We
begin by applying $\Delta_j$ to both sides of the above, using
that it commutes with $S(t)$.  We then insert resolutions of identity
based on the $\Pi$'s and the $\Delta$'s before and after the 
$S_a(s)$ factor, to obtain
\begin{equation}
  \label{eq:ullo}
  \Delta_j u=\Delta_j S(t)u_0-ia\sum_{k,l,m}\int_0^t S(t-s) \Delta_j  \Omega^{-2} \Delta_k 
 \Delta_k\Pi_l S_a(s) \Pi_l \Delta_m \Delta_m f,
\end{equation}

By the endpoint Strichartz estimate for the free group $S(t)$, and generalized H\"older inequality we then have,
 for $n\geq 3$,
\begin{eqnarray*}
\lefteqn{\|\int_0^t S(t-s) J_{jk} M_{kl} S_a(s) N_{lm} \Delta_m f\ ds \|_{L^2_tL^{\frac{2n}{n-2},2}_x}}\qquad\\
&\leq & C \| J_{jk} M_{kl} S_a(s) N_{lm} \Delta_m f \|_{L_t^2L_x^{\frac{2n}{n+2},2}} \\
& \leq & C \|\ |x|^{-1}\|_{L_x^{n,\infty}} \| \ \Omega J_{jk} \Omega \Omega^{-1} M_{kl} \Omega \Omega^{-1}S_a N_{lm} \Delta_m f \|_{L_{t,x}^2}
\end{eqnarray*}
while for $n=2$ we again proceed as before, utilizing Lemma~\ref{lem:CK}, 
Strichartz estimate for the free Schr\"odinger equation with $(p,q)$ as in (\ref{schrodadmis}), and the dual to 
the weighted-$L^2$ estimate~\ref{eq:smoothS} (with $a=0,d=1,\alpha=1/4$) to obtain
\begin{eqnarray*}
\lefteqn{\| S(t) \int S(-s) J_{jk} M_{kl} S_a(s) N_{lm} \Delta_m f ds \|_{L^p_tL^q_x}}\qquad\\
& \leq & C \| \Omega J_{jk} \Omega\Omega^{-1}M_{kl}\Omega\Omega^{-1} S_a(s) N_{lm} 
\Delta_m f \|_{L^2_tL^2_x} 
\end{eqnarray*}
Thus in either case we need to apply Lemma~\ref{xJx} with $\zeta = 1$, Lemma~\ref{x-Mx} with
 $\eta = 1$, Lemma~\ref{lem:Mjk} (all three Lemmas with $d=d_0(n)+1$), and estimate~\ref{eq:smoothS} with
$\alpha=1/4$ to conclude that
\begin{eqnarray}
\|\Delta_j u \|_{L^p_tL^{q,2}_x} & \leq & C \sum_{k,l,m}
2^{-(\epsilon_3|j-k|+\epsilon_2|k-l|+\epsilon_1|l-m|)} \|\Delta_m f\|_{L^2}\nonumber \\
& \leq & C \sum_{m} 2^{-\epsilon|j-m|}\|\Delta_m f\|_{L^2},\label{yow}
\end{eqnarray}
for $\epsilon = \min_{i=1}^3 \epsilon_i < \min\{n/2+d_0(n),\nu_{d_0(n)+1}\}$,
which is the desired frequency-localized version of the endpoint Strichartz estimate, valid for 
data $f \in L^2_{\geq d_0(n)+1}(\RR^n)$.
Here $(p,q) = (2,2n/(n-2))$ if $n\geq 3$ and $(p,q)$ are as in (\ref{schrodadmis}) if $n=2$.

The above estimate can now be used to obtain generalized Strichartz estimates in Sobolev, or more generally, Besov spaces, 
as follows:
Let $(p,q)$ be as in the above and $f\in L^2_{\geq 1+d_0(n)}$. We then have
\begin{eqnarray*} 
\sum_j 2^{2sj}\|\Delta_j u\|_{L^p_t(L^{q,2}_x)}^2 & \leq & C \sum_m 2^{2sm} \|\Delta_m f\|_{L^2}^2 \sum_j 2^{2s(j-m)-2\epsilon|j-m|}\\
& \leq & C \|f\|_{\dot{H}^s}^2 
\end{eqnarray*}
as long as $ |s| < \epsilon$. 

In other words, combining this with Theorem~\ref{thm:strichschrodrad} for $d=d_0(n)+1$ we have proved
\begin{theorem}
  \label{schrodder}
Let  $s$ be such that
\[
\begin{array}{ll} 
-\min\{1+\nu_0,\frac{n}{2},\nu_1\} < s < 
\min\{ 1-\frac{2}{p}+\nu_0,\frac{n}{2},\nu_1\} & \mbox{\rm  if } n\geq 3\\
-\min\{1+\sqrt{a},1+\frac{2}{p}\}<s<\min\{1+\sqrt{a}-\frac{2}{p},1\} & \mbox{\rm if }n=2
\end{array}
\]
where $(p,q)$ are as 
in {\rm(\ref{schrodadmis})}.  Then there exists a constant $C>0$ depending on $n,p,s,a$ 
such that the unique solution $u$ of {\rm(\ref{eq:lse})}
satisfies
\begin{equation}
  \label{eq:strischrodder}
  \|(-\Delta)^{s/2} u\|_{L^p_t(L^q_x)}\leq C
  \|u\|_{L^p_t(\dot{B}^{s,2}_q)}\leq C \|f\|_{\dot{H}^s}.
\end{equation}
\end{theorem}
{\em Proof}: When $n\geq 3$ and $p=2$ this estimate is a simple consequence of the definition of Besov
spaces and (\ref{yow}). Interpolating with the energy estimate
gives the full range of $p$.\cqfd 
\subsection{The wave equation}
We will proceed exactly as in the previous subsection.  We write
\begin{eqnarray}
  \label{eq:ullobis}
  \Delta_j u &=& \dot{W}(t)\Delta_jf+W(t)\Delta_jg\\
&+&a\sum_{k,l,m}\int_0^t W(t-s) \Delta_j  \Omega^{-2} \Delta_k 
 \Delta_k\Pi_l \big[W_a(s) \Pi_l \Delta_m \Delta_m g\nonumber\\
&&\mbox{}+\dot{W}_a(s)\Pi_l\Delta_m\Delta_m f\big]\ ds,\nonumber
\end{eqnarray}
where $W_a(s)$ is the  propagator corresponding to $\partial_t^2+P_a$, i.e. on the $l$-th spherical harmonic subspace
$W_a(s) = \Hi_\nu \frac{\sin (s|\xi|)}{|\xi|} \Hi_\nu$. 

For simplicity, let us assume $f \equiv 0$.
Again, by the endpoint Strichartz estimate \ref{eq:es}, we then have, for $n\geq 4$,
\begin{eqnarray*}
\lefteqn{2^{-\frac{2}{n-1}j}\|\int_0^t W(t-s) J_{jk} M_{kl} W_a(s) N_{lm}
 \Delta_m g\ ds \|_{L^2_tL^{\frac{2(n-1)}{n-3},2}_x}\leq}&&\\
&&C \| J_{jk} M_{kl} W_a(s) N_{lm} \Delta_m g \|_{L_t^2L_x^{\frac{2(n-1)}{n+1},2}} \leq\\
&&C \|\ |x|^{-\frac{n}{n-1}}\|_{L_x^{n-1,\infty}} \| \
\Omega^{\frac{n}{n-1}}J_{jk} \Omega^{\frac{n-2}{n-1}} \ \Omega^{-\frac{n-2}{n-1}}
M_{kl} \Omega^{\frac{n-2}{n-1}}\ \Omega^{-\frac{n-2}{n-1}} W_a N_{lm} \Delta_m g \|_{L_{t,x}^2},
\end{eqnarray*}
Now, we use Lemma \ref{xJx} with $\zeta=\frac{n}{n-1}$, Lemma \ref{x-Mx} with $\eta=\frac{n-2}{n-1}$, Lemma \ref{lem:Mjk}
(all three with $d=d_0(n)+1$), and Corollary~\ref{twabis} with $\alpha=\frac{n-3}{4(n-1)}$ to conclude
\begin{eqnarray}
2^{-\frac{2}{n-1}j}\|\Delta_j u \|_{L^2_tL^{\frac{2(n-1)}{n-3},2}_x} & \leq & C \sum_{k,l,m}
2^{-(\epsilon_3|j-k|+\epsilon_2|k-l|+\epsilon_1|l-m|)}2^{-l\frac{n+1}{2(n-1)}} \|\Delta_m g\|_{L^2}\nonumber \\
& \leq & C \sum_{m} 2^{-\epsilon|j-m|} 2^{-m\frac{n+1}{2(n-1)}}\|\Delta_m g\|_{L^2},\label{yoww}
\end{eqnarray}
for $\epsilon = \min_{i=1}^3 \epsilon_i < \min\{n/2,\nu_1\}-\frac{1}{n-1}$,
which is the desired frequency-localized version of the endpoint Strichartz estimate.

The above can be used to obtain an endpoint Strichartz estimate with derivatives, as follows:  Multiply both 
sides of (\ref{yoww}) by $2^{j(\gamma -\frac{n-3}{2(n-1)})}$, then square both sides, sum over $j$, and change 
the order of summation on the right, to get
\[
\sum_j 2^{2 \sigma j}\|\Delta_j u \|^2_{L^2_tL^{\frac{2(n-1)}{n-3},2}_x} \leq
 C \sum_m 2^{2(\gamma-1)m}\|\Delta_m g\|^2_{L^2}
\]
where $\sigma = \gamma - \frac{n+1}{2(n-1)}$.  This holds provided $|\gamma - \frac{n-3}{2(n-1)}|<\epsilon$.

For the cases $n=2,3$, once again appealing to Lemma~\ref{lem:CK}, we use Strichartz estimate for the free wave equation
to obtain
\[
2^{\sigma j}\|\dot{W}(t)\int W(-s) J_{jk}M_{kl}W_a(s)N_{lm}\Delta_m g ds\|_{L^pL^q}\leq C\|\int W(-s)G(s) ds \|_{\dot{H}^{1/2}},
\]
where $G := J_{jk}M_{kl}W_a(s)N_{lm}\Delta_m g$ and $p,q,\sigma$ are as in (\ref{waveadmis}) with $\gamma = 1/2$. 
Once again, we use the dual to 
(\ref{eq:ff666bis}) with $\alpha = 1/4$ and $a=0$ to obtain
\[
\|\int W(-s)G(s) ds \|_{\dot{H}^{1/2}} \leq C \|\Omega G \|_{L^2(\RR^{n+1})}
\]
Next we use Lemma~\ref{xJx} with $\zeta = 1$, Lemma~\ref{x-Mx} with $\eta = 1$, and (\ref{eq:ff666bis})
with $\alpha=1/4$, all three with $d=d_0(n)+1$ to arrive at
\[
\|\Omega G \|_{L^2(\RR^{n+1})} \leq C 2^{-(\epsilon_3|j-k|+\epsilon_2|k-l|+\epsilon_1|l-m|)}\|P_a^{-1/4}N_{lm}\Delta_m g\|_{L^2}
\]
and from here the proof proceeds as in the case $n\geq 4$ and we obtain
\[
2^{\sigma j} \|\Delta_j u\|_{L^pL^q} \leq C \sum_m 2^{-\epsilon|j-m|} 2^{-m/2} \|\Delta_m g\|_{L^2}
\]
Multiplying by $2^{j(\gamma-\frac{1}{2})}$ and carrying on as above we get
\[
\sum_j 2^{2\sigma' j} \|u\|_{L^pL^q}^2 \leq C \sum_m 2^{2m(\gamma-1)}\|g\|^2_{L^2}
\]
where $\sigma' = \gamma +\frac{1}{p} - n(\frac{1}{2} - \frac{1}{q})$. This holds as long as $|\gamma-\frac{1}{2}|<\epsilon$.

As before, combining the above with Theorem~\ref{thm:strichwaverad} we can deduce generalized
 Strichartz estimates for the wave equation (\ref{eq:lwe}): 
\begin{theorem}
For $n\geq 2$ let  $(p,q,\gamma,\sigma)$
be as in (\ref{waveadmis}). There exists a constant $C>0$ depending on $n,a,p,q,\gamma$ such that
the solution $u$ of {\rm(\ref{eq:lwe})} satisfies
\begin{equation}\label{est:genwav}
\|(-\Delta)^{\sigma/2}u\|_{L^pL^q} \leq C (\|f\|_{\dot{H}^\gamma} + \|g\|_{\dot{H}^{\gamma -1}})
\end{equation}
provided
\[
-\min\{\frac{n-1}{2},\nu_1-\frac{1}{2},1+\nu_0\} < \gamma < 
\min\{ \frac{n+1}{2},\nu_1+\frac{1}{2},1+\nu_0-\frac{1}{p} \}
\]
if $n=2,3$ and
\[
-\min\{\frac{n}{2}-\frac{n+3}{2(n-1)},\nu_1-\frac{n+3}{2(n-1)},1+\nu_0\} < \gamma < 
\min\{ \frac{n+1}{2},\nu_1+\frac{1}{2},1+\nu_0-\frac{1}{p} \}
\]
 if $n\geq 4$.
\end{theorem}

We end this section by giving a nonlinear application for the above estimates.
In \cite{StrTsu97} the authors study the following equation,
\begin{equation}
  \label{eq:nl1}
  \Box u+V(x) u=g(x,u),
\end{equation}
where $V(x)={C}/{|x|^{2- \delta}}$ and $g$ behaves like $|u|^\kappa$ for
 some $\kappa>1$. When $C=0$ (no potential term) the behavior of global
 solutions of small amplitude depend on whether $\kappa$ is larger
 than a critical value $\kappa_c$ (\cite{John79}). In the range $\kappa>\kappa_c$ no
 blow-up occurs. However, adding a potential term can affect the
 solution, leading to blow-up in finite time. Essentially in
 \cite{StrTsu97} blow-up is proved when $\delta>0$ for $\kappa>1$ while global
 existence of smooth solutions is proved when $\delta<0$ for $\kappa>\kappa_c$ and $|C|$ small. 

The inverse-square potential corresponds to the borderline case $\delta= 0$.  
Consider now the following Cauchy problem 
\begin{equation}\label{eq:f1}
\left\{\begin{array}{l}\partial_t^2 u +P_a u=\pm |u|^\kappa,\\
u(0,x) = u_0(x)\\
\partial_t u(0,x) = u_1(x)
\end{array}
\right.
\end{equation}
for $(t,x) \in \mathbb{R}\times\mathbb{R}^n$, where $P_a = -\Delta + \frac{a}{|x|^2}$ with $a> -\lambda(n)^2$
 (See \cite{PST1} for the precise meaning of this operator when $ -\lambda^2 < a < 1-\lambda^2$). 
Set
$\varsigma:=\frac{2}{\kappa-1}$.
 A simple computation shows the equation (\ref{eq:f1}) to
be invariant under the scaling
$u_\varrho=\varrho^{\varsigma}u(\varrho x,\varrho t)$. This
suggests that the equation should be well-posed at the critical level
$\dot{H}^{s_c}$ where $s_c=\frac{n}{2}-\varsigma$, as the
$\dot{H}^{s_c}$ norm of $u$ is left unchanged by this rescaling. Indeed, for
the usual wave equation with no potential, well-posedness holds for
all $s\geq s_c$ when $s_c\geq \frac{1}{2}$, or equivalently $\kappa\geq\frac{n+3}{n-1}$ (
 see \cite{Str81a}, \cite{Kap94} and \cite{LinSog95}).  In \cite{PST1} a similar result was shown to hold
for (\ref{eq:f1}) under the additional assumption that the initial data are {\em radial}. This assumption 
was made only because the linear estimates on which the proof was based were  available only in the radial
case.  The proof of global wellposedness  itself was based on standard iteration arguments and did not use
the assumption of spherical symmetry in any way.  Using the generalized Strichartz estimate (\ref{est:genwav})
that are now available to us we can  remove the assumption of spherical symmetry from this wellposedness 
result.

\begin{theorem}
  \label{ft0}
Let $n\geq 2$,  $\kappa\geq\frac{n+3}{n-1}$, 
$s_c := \frac{n}{2}-\frac{2}{\kappa-1}$.  Suppose $a\in \RR$ is such that 
\[
\sqrt{a+\lambda(n)^2} > \lambda(n) - \frac{2}{\kappa-1} + \max\{\frac{1}{2\kappa},\frac{2}{(n+1)(\kappa-1)}\}.
\]
 Let $(u_0,u_1) \in
(\dot{H}^{s_c},\dot{H}^{s_c-1})$ be  functions with
small norms. Then there exists a unique global  solution
 to (\ref{eq:f1})  such that
 \begin{equation}
   \label{eq:f90}
   u(x,t)\in C_t(\dot{H}^{s_c})\cap L_t^\sigma(\dot{H}^\alpha_q),\qquad\partial_t u(x,t)\in
 C_t(\dot{H}^{s_c-1})\cap
 L_t^\sigma(\dot{H}^{\alpha-1}_q) ,
 \end{equation}
where $q,\sigma,\alpha$ are as follows:
\begin{enumerate}
\item For $n\leq 3$ or $\frac{n+3}{n-1}\leq \kappa \leq \frac{n+1}{n-3}$
\[
\frac{1}{q} = \frac{1}{2} - \frac{4}{(n^2-1)(\kappa-1)},\quad \frac{1}{\sigma} =
\frac{2}{(n+1)(\kappa-1)},\quad \alpha = s_c - \frac{2}{(n-1)(\kappa-1)}.
\]
\item For $n\geq 4$ and $\kappa > \frac{n+1}{n-3}$,
\[
\frac{1}{q} = \frac{1}{2} - \frac{1}{(n-1)\kappa},\quad \frac{1}{\sigma} =
\frac{1}{2\kappa},\quad \alpha = s_c - \frac{n+1}{2(n-1)\kappa}.
\] 
\end{enumerate}
\end{theorem}

{\em Proof:} This result is proved by a contraction mapping argument.  A sequence of Picard iterates is
constructed in the function space
\[
\mathcal{E}= C_t(\dot{H}^{s_c})\cap  L^{\sigma}_{t}(\dot{H}^{\alpha}_{q}).
\]
The nonlinearity in (\ref{eq:f1}) maps $\mathcal{E}$ into another suitably chosen space $\mathcal{F}$, and
the particular choice of the parameters ensures, via the Strichartz estimate (\ref{est:genwav}), 
that $(-\partial_t^2+P_a)^{-1}$ maps
$\mathcal{F}$ back into $\mathcal{E}$.  It is then easy to show that for data of sufficiently small Sobolev
norm, this will be a contraction mapping. See \cite{PST1} for details.

{\bf Acknowledgement:} We would like to thank the referee for several helpful remarks.

\bibliographystyle{plain}

\begin{thebibliography}{10}

\bibitem{BAK}
Matania Ben-Artzi and Sergiu Klainerman.
\newblock Decay and regularity for the {S}chr\"odinger equation.
\newblock {\em J. Anal. Math.}, 58:25--37, 1992.
\newblock Festschrift on the occasion of the 70th birthday of Shmuel Agmon.

\bibitem{CEFG01}
H.~E. Camblong, L.~N. Epele, H.~Fanchiotti, and C.~A. Garcia~Canal.
\newblock Quantum anomaly in molecular physics.
\newblock {\em Phys. Rev. Lett.}, 87(22):220402(4), 2001.

\bibitem{Cas50}
K.~M. Case.
\newblock Singular potentials.
\newblock {\em Physical Rev. (2)}, 80:797--806, 1950.

\bibitem{CheTay82a}
Jeff Cheeger and Michael Taylor.
\newblock On the diffraction of waves by conical singularities. {I}.
\newblock {\em Comm. Pure Appl. Math.}, 35(3):275--331, 1982.

\bibitem{CK}
Michael Christ and Alexander Kiselev.
\newblock Maximal functions associated to filtrations.
\newblock {\em J. Funct. Anal.}, 179(2):409--425, 2001.

\bibitem{GeoVis}
Vladimir Georgiev and Nicola Visciglia.
\newblock Decay estimates for the wave equation with potential.
\newblock {\em preprint}, 2002.

\bibitem{Ginibre}
Jean Ginibre.
\newblock personal communication.

\bibitem{Hos97}
Toshihiko Hoshiro.
\newblock On weighted {$L\sp 2$} estimates of solutions to wave equations.
\newblock {\em J. Anal. Math.}, 72:127--140, 1997.

\bibitem{KJ79}
Arne Jensen and Tosio Kato.
\newblock Spectral properties of {S}chr\"odinger operators and time-decay of
  the wave functions.
\newblock {\em Duke Math. J.}, 46(3):583--611, 1979.

\bibitem{John79}
Fritz John.
\newblock Blow-up of solutions of nonlinear wave equations in three space
  dimensions.
\newblock {\em Manuscripta Math.}, 28(1-3):235--268, 1979.

\bibitem{JouSofSog91}
J.-L. Journ{\'e}, A.~Soffer, and C.~D. Sogge.
\newblock Decay estimates for {S}chr\"odinger operators.
\newblock {\em Comm. Pure Appl. Math.}, 44(5):573--604, 1991.

\bibitem{KalSchWalWus75}
H.~Kalf, U.-W. Schmincke, J.~Walter, and R.~W{\"u}st.
\newblock On the spectral theory of {S}chr\"odinger and {D}irac operators with
  strongly singular potentials.
\newblock In {\em Spectral theory and differential equations (Proc. Sympos.,
  Dundee, 1974; dedicated to Konrad J\"orgens)}, pages 182--226. Lecture Notes
  in Math., Vol. 448. Springer, Berlin, 1975.

\bibitem{Kap94}
Lev Kapitanski.
\newblock Weak and yet weaker solutions of semilinear wave equations.
\newblock {\em Comm. Partial Differential Equations}, 19(9-10):1629--1676,
  1994.

\bibitem{KeeTao98}
Markus Keel and Terence Tao.
\newblock Endpoint {S}trichartz estimates.
\newblock {\em Amer. J. Math.}, 120(5):955--980, 1998.

\bibitem{LinSog95}
Hans Lindblad and Christopher~D. Sogge.
\newblock On existence and scattering with minimal regularity for semilinear
  wave equations.
\newblock {\em J. Funct. Anal.}, 130(2):357--426, 1995.

\bibitem{Mon74}
Vincent Moncrief.
\newblock Odd-parity stability of a {R}eissner-{N}ordstr\"om black hole.
\newblock {\em Phys. Rev. D (3)}, 9:2707, 1974.

\bibitem{Mor68}
Cathleen~S. Morawetz.
\newblock Time decay for the nonlinear {K}lein-{G}ordon equations.
\newblock {\em Proc. Roy. Soc. Ser. A}, 306:291--296, 1968.

\bibitem{One63}
Richard O'Neil.
\newblock Convolution operators and ${L}(p,\,q)$ spaces.
\newblock {\em Duke Math. J.}, 30:129--142, 1963.

\bibitem{PST1}
Fabrice Planchon, John Stalker, and A.~Shadi Tahvildar-Zadeh.
\newblock {$L^p$} estimates for the wave equation with the inverse-square
  potential.
\newblock {\em Discrete Contin. Dynam. Systems}, 9(2):427--442, 2003.

\bibitem{PST2}
Fabrice Planchon, John Stalker, and A.~Shadi Tahvildar-Zadeh.
\newblock Dispersive estimate for the wave equation with the inverse-square
  potential.
\newblock {\em Discrete Contin. Dynam. Systems}, to appear.

\bibitem{ReedSimon3}
Michael Reed and Barry Simon.
\newblock {\em Methods of modern mathematical physics. {III}}.
\newblock Academic Press [Harcourt Brace Jovanovich Publishers], New York,
  1979.
\newblock Scattering theory.

\bibitem{RegWhe57}
Tullio Regge and John~A. Wheeler.
\newblock Stability of a {S}chwarzschild singularity.
\newblock {\em Phys. Rev. (2)}, 108:1063--1069, 1957.

\bibitem{RodS}
Igor Rodnianski and Wilhelm Schlag.
\newblock Time decay for solutions of {S}chr\"odinger equations with rough and
  time dependent potentials.
\newblock {\em preprint}, 2001.

\bibitem{Sim92}
Barry Simon.
\newblock Best constants in some operator smoothness estimates.
\newblock {\em J. Funct. Anal.}, 107(1):66--71, 1992.

\bibitem{Str81a}
Walter~A. Strauss.
\newblock Nonlinear scattering theory at low energy.
\newblock {\em J. Funct. Anal.}, 41(1):110--133, 1981.

\bibitem{StrTsu97}
Walter~A. Strauss and Kimitoshi Tsutaya.
\newblock Existence and blow up of small amplitude nonlinear waves with a
  negative potential.
\newblock {\em Discrete Contin. Dynam. Systems}, 3(2):175--188, 1997.

\bibitem{Tit46}
E.~C. Titchmarsh.
\newblock {\em Eigenfunction Expansions Associated with Second-Order
  Differential Equations}.
\newblock University Press, Oxford, 1946.

\bibitem{VazZua00}
Juan~Luis Vazquez and Enrike Zuazua.
\newblock The {H}ardy inequality and the asymptotic behaviour of the heat
  equation with an inverse-square potential.
\newblock {\em J. Funct. Anal.}, 173(1):103--153, 2000.

\bibitem{Yaj95b}
Kenji Yajima.
\newblock The ${W}\sp {k,p}$-continuity of wave operators for {S}chr\"odinger
  operators.
\newblock {\em J. Math. Soc. Japan}, 47(3):551--581, 1995.

\bibitem{Zer70}
Frank~J. Zerilli.
\newblock Gravitational field of a particle falling in a {S}chwarzschild
  geometry analyzed in tensor harmonics.
\newblock {\em Phys. Rev. D (3)}, 2:2141--2160, 1970.

\bibitem{Zer74}
Frank~J. Zerilli.
\newblock Perturbation analysis for gravitational and electromagnetic radiation
  in a {R}eissner-{N}ordstr\"om geometry.
\newblock {\em Phys. Rev. D (3)}, 9:860, 1974.

\end{thebibliography}
\def\cprime{$'$}

\vskip 1cm 
\noindent { $^{\mathbf   a}$ 
 D\'epartement de Math{\'e}matiques\\
Batiment 425, Universit\'e Paris-Sud\\
F-91405 Orsay Cedex\\
 Email: {\sl Nicolas.Burq@math.u-psud.fr}
\ \\
 \\ $^{\mathbf b}$
Laboratoire Analyse, G\'eom\'etrie \& Applications\\
UMR 7539, Institut Galil\'ee\\
Universit\'e Paris 13\\
99 avenue J.B. Cl\'ement\\
F-93430 Villetaneuse\\
 Email: {\sl Fab@math.univ-paris13.fr}
\ \\
 \\ $^{\mathbf c}$
Department of Mathematics\\
Princeton University,\\
 Princeton NJ 08544\\
Email: {\sl Stalker@math.princeton.edu}
\ \\
 \\ $^{\mathbf d}$
Department of Mathematics\\
Rutgers, The State University of New Jersey\\
 110 Frelinghuysen Road, Piscataway NJ 08854\\
Email: {\sl Shadi@math.rutgers.edu}
 }
\end{document}